\documentclass[journal]{IEEEtran}
\usepackage{amsmath,amssymb,amsfonts}
\usepackage{algorithmic}
\usepackage{graphicx}
\usepackage{subcaption}
\usepackage{algorithm,algorithmic}
\usepackage{hyperref}
\usepackage{amsthm}
\usepackage{subcaption}
\usepackage[table]{xcolor}
\usepackage{chngcntr}
\usepackage{booktabs} 
\usepackage{multirow}
\usepackage{threeparttable}
\newtheorem{theorem}{Theorem}%{Theorem}%[section]
\newtheorem{lemma}{Lemma}
\newtheorem{corollary}{Corollary}
\newtheorem{assumption}{Assumption}
\newtheorem{remark}{Remark}
\usepackage{pifont}
\usepackage{booktabs}
\newcommand{\cmark}{\ding{51}}
\newcommand{\xmark}{\ding{55}}

\usepackage{xcolor}

\begin{document}

\title{Variance Reduction for Decentralized Stochastic Optimization with Multiple Local Updates
}

\author{Jie Liu, Nicola Bastianello, \IEEEmembership{Member, IEEE}, Wenying Xu, \IEEEmembership{Senior Member, IEEE}, Daniel W. C. Ho, \IEEEmembership{Life Fellow, IEEE}, and Karl H. Johansson, \IEEEmembership{Fellow, IEEE}
	 \thanks{%This work was supported in part by the National Natural Science Foundation of China under Grant 62603151, 62573121. The work of K. H. Johansson was funded by the Swedish Research Council Distinguished Professor Grant 2017-01078 Knut and Alice Wallenberg Foundation Wallenberg Scholar Grant. (Corresponding author: Wenying Xu).
     
     J. Liu and W. Xu  are  with the School of Mathematics, Southeast University, Nanjing 211189, China (jliu@seu.edu.cn; wyxu@seu.edu.cn). N. Bastianello is with the Department of Electronic and Electrical Engineering, University College London, United Kingdom (n.bastianello@ucl.ac.uk). Daniel W. C. Ho is  with the Department of Mathematics, City University of Hong Kong, Hong Kong, China (madaniel@cityu.edu.hk). K. H. Johansson is with the School of Electrical Engineering and Computer Science, and Digital Futures, KTH Royal Institute of Technology, 11428 Stockholm, Sweden (kallej@kth.se).}
}

%\author{wait \thanks{wait}    }

\markboth{ }%
{Shell \MakeLowercase{\textit{et al.}}: Bare Demo of IEEEtran.cls for IEEE Journals}

\maketitle

\begin{abstract}
Multiple local updates are widely used in decentralized optimization to enhance communication efficiency, often in conjunction with stochastic gradients when exact gradient computation is costly. However, such algorithms generally converge only to a neighborhood of the optimum, whose size depends on the stepsize and the variance of gradient noise. Variance reduction has been incorporated into ADMM-based methods to achieve exact convergence in expectation, but their stepsizes typically scale inversely with both the square root of the network size and the cube of the number of local updates, leading to slow convergence as the network size or the number of local updates increases. In this paper, we propose a novel primal-dual decentralized optimization algorithm that integrates variance reduction with multiple local updates, with a stepsize independent of the network size and inversely proportional to the number of local updates (rather than its cube), thereby enabling significantly larger stepsizes and faster convergence. Moreover, to the best of our knowledge, we establish the first linear convergence guarantee in expectation for decentralized optimization with multiple local updates under stochastic gradients in the strongly convex setting. Beyond the theoretical guarantee, the proposed method exchanges only one variable per round, halving the communication cost compared with gradient-tracking methods. Numerical experiments demonstrate the advantages of the proposed algorithm in terms of convergence speed and accuracy.

%In addition, the proposed method exchanges only one variable per round, halving the per-round communication cost compared with gradient-tracking-based methods that exchange two.  Numerical experiments on distributed estimation problems confirm the superior convergence speed and accuracy of the proposed method over existing approaches.

%Variance reduction provides a natural approach to mitigating gradient noise and enabling exact convergence in expectation, and has been incorporated into existing ADMM-based methods. Nevertheless, these methods typically require stepsizes that scale inversely with both the square root of the network size and the cube of the number of local updates. Consequently, the stepsize becomes prohibitively small as either quantity increases, leading to slow convergence.

%In addition, each agent communicates only a single variable per round, thereby halving the per-round communication cost compared with gradient-tracking-based methods that exchange both optimization and auxiliary variables.
\end{abstract}


\IEEEpeerreviewmaketitle

\section{Introduction}
Decentralized optimization has emerged as a compelling paradigm with broad applications in artificial intelligence, including  neural network training \cite{pmlr-v54-mcmahan17a}, regression analysis \cite{yuan2020distributed}, and distributed estimation \cite{li2023generalized}. In this setting, each agent performs local computation and coordinates exclusively with its neighboring agents through information exchange. Relative to centralized architectures, decentralized optimization offers several distinctive advantages, including enhanced robustness against single-point failures, improved preservation of data privacy, and superior scalability and communication efficiency \cite{yao_yang_survey}. These merits have motivated extensive research on decentralized optimization algorithms, with particular emphasis on their theoretical convergence guarantees \cite{Nedich_subgradient,xiuxianli1,xin2020variance1,Leishi1,ED_Yuan,alghunaim2022unified,liu2022accelerated,bo2024quantization}.

Communication among agents plays a fundamental role in decentralized optimization. However, in large networks, the associated communication overhead can become a major bottleneck. Consequently, various communication-efficient techniques, including compression and quantization, have been developed for decentralized optimization \cite{Leishi1,liu2024distributed,11353935,11165038}. Notably, a highly effective approach is to let each agent perform multiple local updates between communication rounds, a strategy that has achieved remarkable success in federated learning by significantly reducing communication rounds without compromising convergence \cite{pmlr-v54-mcmahan17a,A_Mitra_Linear,pmlr-v119-karimireddy20a,liu2026achieving}. The success of this strategy in federated learning naturally raises the following question: can it be extended to fully decentralized optimization? This extension, however, is far from straightforward. Without a central server to coordinate updates, agents in the fully decentralized setting must reach consensus through iterative information exchange with their neighbors using a mixing matrix. During the local update phase, however, they rely solely on local gradient information, which introduces significant challenges in both algorithm design and theoretical analysis. Several existing works have established convergence guarantees for multiple local updates in decentralized optimization under exact gradients \cite{liu2025guaranteeing,Mind_jieliu,LOCAL_SMALL_STEPSIZE,kunyuan1,wu2025effectiveness}. Furthermore, via the Performance Estimation Problem (PEP) framework it has been shown that for the the seminal DIGing algorithm multiple local updates can indeed accelerate convergence, thereby enhancing communication efficiency \cite{wang2026provable}.

In modern large-scale applications, however, computing the exact gradient (full-batch) is computationally prohibitive, making stochastic gradients (typically via mini-batch sampling) the standard choice 
\cite{pmlr-v54-mcmahan17a,Shi_Pu_Tracking,xin2020variance1,xin2020variance2}. Although decentralized stochastic gradient methods have achieved remarkable empirical success, the inherent noise in stochastic gradients prevents exact convergence. The convergence is only to a neighborhood of the optimum, whose radius is determined by the stepsize $\alpha$ and the gradient variance $\sigma^2$ \cite{liu2025guaranteeing,Mind_jieliu,LED,taolin1,LOCAL1_AVERAGING,huang2025beyond,huang2025distributed,kunyuan_local,shipu_local_update,time_varying_local}. To address this issue, variance reduction techniques have been extended to decentralized settings to achieve exact convergence in expectation, where the variance of the gradient estimators is bounded by consensus errors rather than a constant \cite{xin2020variance1,xin2020variance2}. With appropriately chosen stepsizes, these algorithms thus achieve exact convergence in expectation \cite{xin2020variance1,xin2020variance2}. However, these results rely on a single local update per communication round, and their extension to the multiple-local-update setting remains largely unexplored.

Decentralized optimization with both stochastic gradients and multiple local updates converges only to a neighborhood of size $\mathcal{O}(\alpha\sigma^2)$, where $\alpha$ is the stepsize and $\sigma^2$ is the stochastic gradient variance \cite{liu2025guaranteeing,Mind_jieliu,LED,taolin1}. Incorporating variance reduction into this setting, however, poses substantial challenges. The fundamental difficulty is that during the multiple local update phase, each agent updates its local variable using only its local stochastic gradient, rather than the global gradient, so the heterogeneity of local objective functions hinders consensus during this process \cite{pmlr-v119-karimireddy20a,A_Mitra_Linear,A_Mitra_Aggregated}. Applying variance reduction in decentralized settings typically requires bounding the variance of stochastic gradients by consensus errors \cite{xin2020variance1,xin2020variance2}. This requirement, combined with the difficulty of maintaining consensus across multiple local updates, inevitably forces a sufficiently small stepsize to guarantee convergence \cite{ren2025communication}. Specifically, the stepsize in \cite{ren2025communication} scales as $1/(\sqrt{N}\tau^3)$, thus inversely proportional to both the square root of the network size $N$ and the cube of the number of local updates $\tau$. For large networks or when $\tau$ is large, this stepsize becomes small, leading to slow convergence.

To address the above limitations, this paper proposes a primal-dual decentralized algorithm that integrates variance reduction with multiple local updates. The main contributions are summarized as follows:

\begin{itemize}
    \item \textbf{New Algorithm Design.} We propose \textbf{LPDVR} (\textbf{L}ocal \textbf{P}rimal-\textbf{D}ual-based Algorithm with \textbf{V}ariance \textbf{R}eduction), a novel decentralized algorithm that integrates variance reduction with multiple local updates. The key idea is to construct a variance-reduced stochastic gradient estimator upon the local primal variable maintained during local updates, rather than the optimization variable. This design ensures that the gradient variance remains controlled across multiple local steps and effectively counteracts the heterogeneity-induced drift during local updates. As a result, LPDVR achieves exact convergence in expectation for nonconvex, convex, and strongly convex costs, while also enabling large stepsizes and high communication efficiency. Section~\ref{comparisons_with_many_existing_works} provides a detailed comparison of convergence properties, stepsize design, and communication overhead among LPDVR and representative existing methods \cite{liu2025guaranteeing,Mind_jieliu,LOCAL_SMALL_STEPSIZE,kunyuan1,wu2025effectiveness,LED,taolin1,huang2025beyond,huang2025distributed,ren2025communication,WeiShi3}. Numerical experiments confirm that LPDVR attains faster convergence and higher accuracy than existing methods.
    
    \item \textbf{Enhanced Convergence Guarantees.} For most existing stochastic gradient methods with multiple local updates \cite{LED,taolin1,huang2025beyond,liu2025guaranteeing,Mind_jieliu}, the established convergence guarantees are limited. With stepsize $\alpha$ and stochastic gradient variance $\sigma^2$, both nonconvex and convex objectives achieve an $\mathcal{O}(1/T)$ rate to a neighborhood of radius $\mathcal{O}(\alpha\sigma^2)$ \cite{LED,taolin1,huang2025beyond,huang2025distributed}. The same limitation applies to strongly convex cases and nonconvex cases satisfying the Polyak--Lojasiewicz condition, where existing methods only achieve linear convergence to a neighborhood rather than the exact optimum \cite{LED,huang2025beyond,huang2025distributed}. In sharp contrast, LPDVR achieves exact convergence at a $\mathcal{O}(1/T)$ rate in expectation for both the nonconvex and convex cases. Furthermore, even compared with the recent variance-reduced method \cite{ren2025communication} that achieves exact convergence only in nonconvex and convex cases, LPDVR attains exact linear convergence for strongly convex objectives, a significant theoretical advance that, to the best of our knowledge, has not been established in prior work.

    \item \textbf{Improved Stepsize Design.} The stepsize of LPDVR scales as $1/\tau$ and is independent of the network size. Notably, this stepsize condition is comparable to, or even less restrictive than, those required by conventional stochastic or exact gradient methods \cite{liu2025guaranteeing,LOCAL_SMALL_STEPSIZE,kunyuan1,LED,taolin1,WeiShi3,huang2025beyond,huang2025distributed}, confirming that variance reduction is achieved without enforcing a more conservative stepsize. This is made possible by a novel variance reduction mechanism built upon the local primal variable, combined with refined mathematical derivations that tightly bound the accumulated errors during multiple local updates. In sharp contrast, the recent variance-reduced ADMM-based method \cite{ren2025communication} requires a stepsize scaling as $1/(\sqrt{N}\tau^3)$. Our method offers two significant advantages: (i) it removes the dependence on $N$, making it suitable for large-scale networks; (ii) it substantially improves the dependence on $\tau$ from $1/\tau^3$ to $1/\tau$, enlarging the feasible stepsize especially when a large number of local updates is performed.

    \item \textbf{Efficient Communication.} Based on the design of LPDVR, each agent shares only its primal variable with neighbors at each communication round, without transmitting any auxiliary variables. This is in sharp contrast to many existing decentralized algorithms \cite{LOCAL_SMALL_STEPSIZE,kunyuan1,taolin1,LOCAL1_AVERAGING,WeiShi3,Nedich_time_varying2}, which require the exchange of two variables per round, typically the optimization variable and an additional auxiliary variable. By cutting the per-round communication overhead in half, LPDVR achieves significantly higher communication efficiency, making it particularly suitable for bandwidth-limited environments and high-dimensional optimization and learning applications.
\end{itemize}

The remainder of this paper is organized as follows.
Section~\ref{problem_section} introduces the problem formulation and basic assumptions.
Section~\ref{step_details} presents the proposed LPDVR algorithm and its matrix form.
Section~\ref{Preliminaries_for_CA} establishes preliminary results for the convergence analysis.
Section~\ref{convergence_for_nonconvex_case} analyzes the convergence of LPDVR in the nonconvex setting, while
Section~\ref{convergence_analysis_for_c_and_sc} considers the convex and strongly convex settings.
Section~\ref{comparisons_with_many_existing_works} compares LPDVR with existing methods, highlighting its advantages in stepsize requirements and convergence guarantees.
Section~\ref{numerical_evaluations} presents numerical experiments to demonstrate the performance advantages of LPDVR over existing methods.
Finally, Section~\ref{conclusion} concludes the paper.
Additional proofs and supporting lemmas are provided in the appendices.

\textbf{Notations:} Throughout this paper, \(\mathbb{Z}^{+}\), \(\mathbb{R}^n\), and \(\mathbb{R}^{m \times n}\) denote the sets of positive integers, real \(n\)-dimensional vectors, and real \(m \times n\) matrices, respectively. The symbols \(\mathbf{1}_N \in \mathbb{R}^N\) and \(\mathbf{I}_N \in \mathbb{R}^{N \times N}\) denote the all-one vector and the identity matrix, respectively. For a vector \(x \in \mathbb{R}^n\) and a matrix \(A \in \mathbb{R}^{N \times N}\), we denote the \(j\)-th entry of \(x\) and the \((i,j)\)-th entry of \(A\) by \([x]_j\) and \([A]_{ij}\), respectively. For \(x \in \mathbb{R}^n\), \(\|x\|\) denotes its Euclidean norm \(\sqrt{\sum_{j=1}^n [x]_j^2}\). For a matrix \(A \in \mathbb{R}^{N \times N}\), \(\|A\|\) denotes the norm of \(A\) induced by the Euclidean norm, \(\|A\| = \sup_{\|x\|=1} \|Ax\|\). Additionally, the induced matrix \(1\)-norm of \(A\) is given by \(\|A\|_1 = \max_{1 \leq j \leq N} \sum_{i=1}^N |[A]_{ij}|\). Finally, \(\operatorname{diag}\{a_i\}_{i=1}^{N}\) denotes the \(N \times N\) diagonal matrix with diagonal entries \(a_1, \ldots, a_N\) and zeros elsewhere.

\section{Problem Setting}\label{problem_section}
We consider an undirected network \(\mathcal{G} = \{\mathcal{V}, \mathcal{E},W\}\), where \(\mathcal{V} = \{1, \dots, N\}\) is the set of agents and \(\mathcal{E} \subseteq \mathcal{V} \times \mathcal{V}\) is the set of edges. An edge \((i,j) \in \mathcal{E}\) indicates that agents \(i\) and \(j\) can exchange information with each other. We investigate the following decentralized finite-sum minimization
problems over \( \mathcal{G} = \{\mathcal{V}, \mathcal{E},W\} \):
\begin{align}\label{op_problem}
\min_{x \in \mathbb{R}^n} f(x) = \frac{1}{N} \sum_{i=1}^{N} f_i(x),  \quad f_i(x)=\frac{1}{m_i}\sum^{m_i}_{j=1}f_{i,j}(x),
\end{align}
where $f_i: \mathbb{R}^n \rightarrow \mathbb{R}$ is the local objective function private to agent $i \in \mathcal{V}$, with each component $f_{i,j}$ corresponding to a single data point. We impose the following assumption on the component functions:
\begin{assumption}\label{smooth_assumption}
The component function $f_{i,j}(x)$ is $L$-smooth over $\mathbb{R}^n$, i.e., $\Vert \nabla f_{i,j}(x)-\nabla f_{i,j}(y)\Vert\leq L \Vert x-y\Vert$ for any $x,y\in\mathbb{R}^n$, $i\in\mathcal{V}$, and $1\leq j\leq m_i$.
\end{assumption}
Under Assumption~\ref{smooth_assumption} and the definitions in \eqref{op_problem}, it follows that for any $i\in\mathcal{V}$, both the global objective $f(x)$ and the local objective function $f_i(x)$ are $L$-smooth over $\mathbb{R}^n$. Additionally, we assume the following:
\begin{assumption}\label{non_empty}
The optimal solution set $\mathcal{X}^*=\{x^*\in\mathbb{R}^n |x^*=\arg\min_{x\in\mathbb{R}^n} f(x)\}$ is nonempty. 
\end{assumption}
Assumptions \ref{smooth_assumption} and \ref{non_empty} are standard in decentralized nonconvex optimization \cite{LED,taolin1,Ming_Yan2,xin2020variance2}. The graph $\mathcal{G}$ is assumed to be undirected and connected, and its mixing matrix $W$ is nonnegative and satisfies $[W]_{ij}=0$ whenever $i\neq j$ and $(i,j)\notin\mathcal{E}$. We also impose the following standard assumption on the communication graph \cite{alghunaim2022unified,LED,Ming_Yan2,WeiShi1,Ming_Yan1}.
\begin{assumption}\label{mixing_matrix} 
The mixing matrix $W$ is positive definite, symmetric, doubly stochastic, and primitive.
\end{assumption}
It is worth noting that the positive definiteness in Assumption \ref{mixing_matrix} is easy to satisfy. Indeed, for any symmetric, doubly stochastic, and primitive matrix $W$, the matrix $\frac{\mathbf{I}_N+W}{2}$ satisfies Assumption \ref{mixing_matrix}. Under Assumption \ref{mixing_matrix}, the eigenvalues of $W$ are real and satisfy $1=\lambda_1>\lambda_2\geq \cdots \geq \lambda_N>0$, and the network mixing rate is given by
\begin{align}\label{eigenvalues_property}
\lambda=\Vert W-\frac{1}{N}\mathbf{1}_N\mathbf{1}^{\top}_N\Vert<1,
\end{align} 
which coincides with $\lambda_2$. The mixing matrix $W$ critically affects the convergence behavior of the proposed algorithm. We summarize several key properties of $W$ below.
\begin{lemma}[\cite{LED}]\label{lemma_from_LED1}
Under Assumption \ref{mixing_matrix}, there exist $\Lambda=\operatorname{diag}\{\lambda_i\}^N_{i=2}$, $Q\in\mathbb{R}^{N\times (N-1)}$ such that $QQ^{\top}=\mathbf{I}_N-\frac{1}{N}\mathbf{1}_N\mathbf{1}^{\top}_N$ and $Q^{\top}Q=\mathbf{I}_{N-1}$, and consequently, we have $\mathbf{W}=\frac{1}{N}\mathbf{1}_N\mathbf{1}^{\top}_N\otimes \mathbf{I}_n+\mathbf{Q}\mathbf{\Lambda}\mathbf{Q}^{\top}$ and $\mathbf{Q}^{\top}\mathbf{W}=\mathbf{\Lambda}\mathbf{Q}^{\top}$, where $\mathbf{W}=W\otimes \mathbf{I}_n$, $\mathbf{\Lambda}=\Lambda\otimes \mathbf{I}_n$, and $\mathbf{Q}=Q\otimes \mathbf{I}_n$. Moreover, we have the identities $\Vert \mathbf{Q}\Vert=1$ and $\Vert (\mathbf{I}_{n(N-1)}-\mathbf{\Lambda})^{-1}\Vert=\frac{1}{1-\lambda}$.
\end{lemma}
Next, we present our algorithm and characterize its convergence.

\section{Algorithm Description}\label{step_details}
In this section, we introduce \textbf{LPDVR} (\textbf{L}ocal \textbf{P}rimal-\textbf{D}ual-based Algorithm with \textbf{V}ariance \textbf{R}eduction), a decentralized primal-dual algorithm that seamlessly combines variance reduction with multiple local updates, achieving exact convergence in expectation without sacrificing communication efficiency.

\subsection{Notation}

We start by defining the notation used in the algorithm. Each agent $i \in \mathcal{V}$ maintains an optimization variable $x_i(t) \in \mathbb{R}^n$ at iteration $t$, initialized as $x_i(0) = \mathbf{0}_n$, and a dual variable $y_i(t) \in \mathbb{R}^n$, initialized as $y_i(0) = \mathbf{0}_n$. To facilitate variance reduction, each agent stores $z_{i,j}(t)$, which denotes the latest iterate at which the local component gradient $\nabla f_{i,j}$ was evaluated up to iteration $t$, and computes a gradient estimator based on these stored values. During the local update phase, the local primal variable is denoted by $\phi_i(t, r)$, where $r$ indexes the local update steps within iteration $t$. The parameters $\tau\in \mathbb{Z}^{+}$ and $\alpha>0$ denote the number of local updates per round and the stepsize, respectively. We now detail the proposed algorithm. 

\subsection{Algorithm Details}

\begin{algorithm}[t]
\caption{LPDVR}
\label{algorithm_recur}
\begin{algorithmic}
\STATE \textbf{Initialization:} $x_i(0)=y_i(0)=\mathbf{0}_n$, $\alpha>0$, $\tau$, and  $\beta=\frac{1}{\tau}$.
\FOR{$t=0$ \textbf{to} $T$}
    \FOR{each agent $i=1,2,\ldots,N$ in parallel}
        \STATE Set $\phi_i(t,0)=x_i(t)$ and $z_{i,h}(t,0)=x_i(t)$ for all $h$.
        \FOR{$r=0$ \textbf{to} $\tau-1$}
            \STATE Sample $s_i(t,r)$ uniformly from $\{1,\dots,m_i\}$, compute the gradient estimator $g_i(\phi_i(t,r))$ via \eqref{gradient_process_algorithm1}, and perform the local primal update:
            \begin{align}
            \hspace{-0.19cm}\phi_i(t,r+1) = \phi_i(t,r) - \alpha g_i(\phi_i(t,r)) - \beta y_i(t).\label{update_phi}
            \end{align}
            Update the selected stored iterate:
            \begin{align}\label{add_for_if_algorithm}
                z_{i,h}(t,r+1) = 
                \begin{cases}
                    \phi_i(t,r+1), & h = s_i(t,r), \\
                    z_{i,h}(t,r), & h \neq s_i(t,r).
                \end{cases}
            \end{align}
        \ENDFOR
        \STATE Update variables based on neighbors' information:
        \vspace{-0.2cm}
        \begin{align}
            x_i(t+1) &= \sum_{j=1}^{N} [W]_{ij} \phi_j(t,\tau),\label{update_y} \\
            y_i(t+1) &= y_i(t) + \phi_i(t,\tau) - x_i(t+1).\label{update_x}
        \end{align}
    \ENDFOR
\ENDFOR
\end{algorithmic}
\end{algorithm}
The proposed LPDVR algorithm is summarized in Algorithm~\ref{algorithm_recur}. At the beginning of each outer iteration $t$, agent $i$ initializes $\phi_i(t,0)=x_i(t)$ and resets $z_{i,h}(t,0)=x_i(t)$ for all $h$. In the local update, for $r=0,\dots,\tau-1$, each agent samples $s_i(t,r)$ uniformly from $\{1,\dots,m_i\}$, computes the variance-reduced gradient estimator
\begin{align}
g_i(\phi_i(t,r)) = &\frac{1}{m_i}\sum_{h=1}^{m_i}\nabla f_{i,h}(z_{i,h}(t,r))+\nabla f_{i,s_i(t,r)}(\phi_i(t,r))\nonumber \\
& - \nabla f_{i,s_i(t,r)}(z_{i,s_i(t,r)}(t,r)), \label{gradient_process_algorithm1}
\end{align}
and updates the local primal variable $\phi_i(t,r+1)$ according to~\eqref{update_phi}. The selected stored iterate $z_{i,s_i(t,r)}$ is then updated to $\phi_i(t,r+1)$, while the others remain unchanged.

After $\tau$ local updates, each agent exchanges $\phi_i(t,\tau)$ with its neighbors. The optimization variable $x_i(t+1)$ is updated as a weighted average of the neighbors' local variables, and the dual variable $y_i(t+1)$ is updated to accumulate the consensus discrepancy, as shown in~\eqref{update_y} and~\eqref{update_x}.

\begin{remark}
LPDVR can be regarded as a primal-dual algorithm for decentralized consensus optimization. Each agent $i$ simultaneously maintains a primal variable $\phi_i(t,r)$, representing the local decision iterate, and a dual variable $y_i(t)$, acting as the Lagrange multiplier associated with the consensus constraint. Within each communication round, the algorithm alternates two steps: (i) the primal update~\eqref{update_phi} performs stochastic gradient descent on the local objective augmented with the dual-feedback term $-\beta y_i(t)$; and (ii) the dual update~\eqref{update_x} executes a multiplier ascent on the consensus violation $\phi_i(t,\tau) - x_i(t+1)$. This \emph{primal-descent and dual-ascent} alternation is the defining feature of primal-dual algorithms. The structure of LPDVR is inspired by \cite{LED}. However, unlike \cite{LED}, LPDVR achieves exact convergence in expectation under stochastic gradients without enforcing a more restrictive stepsize, whereas the method in \cite{LED} converges only to a neighborhood of size $\mathcal{O}(\alpha\sigma^2)$, where $\alpha$ is the stepsize and $\sigma^2$ is the variance of the stochastic gradient. A more detailed comparison is provided in Section~\ref{comparisons_with_many_existing_works}.
\end{remark}

\subsection{Matrix Formulation}
For the convergence analysis of LPDVR, we first present its matrix form. For notational convenience, for any vector $X=\mathrm{col}\{x_i\}_{i=1}^{N}\in\mathbb{R}^{nN}$, we define
\begin{align}
\hspace{-0.17cm}\nabla F(X) \triangleq \mathrm{col}\{\nabla f_i(x_i)\}_{i=1}^{N},\ \overline{\nabla F}(X)=\frac{1}{N}\sum^{N}_{i=1}\nabla f_i(x_i).
\end{align}
Moreover, we introduce the following stacked vectors of the iterates and gradient estimates:
\begin{equation}\label{matrix_definition}
\hspace{-0.24cm}\left\{
\begin{aligned}
&X(t) = \mathrm{col}\{x_i(t)\}_{i=1}^{N}, \Phi(t,r)= \mathrm{col}\{\phi_i(t,r)\}_{i=1}^{N},\\
&Y(t) = \mathrm{col}\{y_i(t)\}_{i=1}^{N}, G(\Phi(t,r)) = \mathrm{col}\{g_i(\phi_i(t,r))\}_{i=1}^{N}.\hspace{-0.2cm}
\end{aligned}
\right.
\end{equation}
Rewriting \eqref{update_phi} in matrix form gives
\begin{align}\label{phi_matrix_from}
\Phi(t,\tau)=X(t)-\alpha\sum^{\tau-1}_{l=0}G(\Phi(t,l))- Y(t),
\end{align}
where we have used the identity $\beta\tau = 1$ to simplify the expression. Combining \eqref{update_phi}, \eqref{update_y}, \eqref{update_x}, \eqref{matrix_definition} with \eqref{phi_matrix_from} yields the following equivalent form of LPDVR:
\begin{align}
X(t+1)=&\mathbf{W}\Phi(t,\tau),\label{matrix_form_algorithm_x}\\
Y(t+1)=&Y(t)+(\mathbf{I}-\mathbf{W})\Phi(t,\tau).\label{matrix_form_algorithm_y}
\end{align}
We further define the average of the optimization variable
\begin{align}\label{averaging_x}
\bar{x}(t)=\frac{1}{N}\sum^{N}_{i=1}x_i(t),\quad \bar{X}(t) = \mathbf{1}_N \otimes \bar{x}(t)
\end{align}
Based on \eqref{averaging_x}, we introduce the auxiliary variable
\begin{align}\label{update_rule_E}
E(t)=\mathrm{col}\{e_i(t)\}_{i=1}^{N} \triangleq Y(t)+\alpha\tau\nabla F(\bar{X}(t))
\end{align}
to facilitate the consensus error analysis. Recall from the update rule of $Y(t)$ in \eqref{matrix_form_algorithm_y} and the initialization $Y(0)=\mathbf{0}_{nN}$ that $(\mathbf{1}^{\top} \otimes I_n)Y(t)=0$ holds. Averaging \eqref{update_rule_E} over all agents gives
\begin{align}
\bar{e}(t) \triangleq \frac{1}{N}\sum_{i=1}^{N} e_i(t) = \frac{\alpha\tau}{N}\sum_{i=1}^{N} \nabla f_i(\bar{x}(t)). \label{average_e}
\end{align}
This relation indicates that the average of $e_i(t)$ effectively tracks the average gradient across agents, motivating us to analyze the convergence of LPDVR through the properties of $E(t)$. To this end, we define the auxiliary variable
\begin{align}\label{matrix_definition2}
\hspace{-0.4em}S(t,r)=\mathrm{col}\{s_i(t,r)\}^{N}_{i=1} \triangleq G(\Phi(t,r))-\nabla F(\Phi(t,r)).
\end{align}
Then, the update rule for $Y(t)$ in \eqref{matrix_form_algorithm_y} can be reformulated as
\begin{align}
E(t&+1)=\mathbf{W}E(t)+\alpha\tau\big(\nabla F(\bar{X}(t+1))-\nabla F(\bar{X}(t))\big)\nonumber\\
&+(\mathbf{I}-\mathbf{W})\big\{X(t)-\alpha\sum^{\tau-1}_{r=0}\big(\nabla F(\Phi(t,r))-\nabla F(\bar{X}(t))\nonumber\\
&+S(t,r)\big)\big\}.\label{E_update_rule}
\end{align}
Applying the decomposition of $\mathbf{W}$ in Lemma~\ref{lemma_from_LED1} to \eqref{matrix_form_algorithm_x} and \eqref{E_update_rule} yields the following lemma.
\begin{lemma}\label{lemma_most_important}
The dynamics of $X(t)$ and $E(t)$ in \eqref{matrix_form_algorithm_x} and \eqref{E_update_rule} admit the following equivalent form:
\begin{align}\label{transformed_system}
D(t+1)= \Delta D(t) - \alpha \bigl(S(t) + H(t)\bigr),
\end{align}
where $\Delta$ is a matrix with norm $\delta \triangleq \|\Delta\| = \sqrt{\lambda} < 1$, 
\begin{align*}
\mathbf{V} &=
\begin{bmatrix}
\mathbf{\Lambda}^{\frac{1}{2}} & \mathbf{\Lambda}^{\frac{1}{2}} \\
-\iota \mathbf{B}^{\frac{1}{2}} & \iota \mathbf{B}^{\frac{1}{2}}
\end{bmatrix} \mathbf{P}^{\top},\ 
D(t) = \mathbf{V}^{-1}
\begin{bmatrix}
\mathbf{Q}^{\top} \mathbf{X}(t) \\
\mathbf{Q}^{\top} \mathbf{E}(t)
\end{bmatrix},\\
S(t) &= \mathbf{V}^{-1}
\begin{bmatrix}
\mathbf{\Lambda} \mathbf{Q}^{\top} \sum_{r=0}^{\tau-1}S(t,r) \\
\mathbf{B} \mathbf{Q}^{\top} \sum_{r=0}^{\tau-1} S(t,r)
\end{bmatrix}, H(t) = \mathbf{V}^{-1}
\begin{bmatrix}
H_1(t) \\
H_2(t)
\end{bmatrix},
\end{align*}
  $H_1(t)=\mathbf{\Lambda} \mathbf{Q}^{\top} \sum_{r=0}^{\tau-1} \bigl(\nabla F(\Phi(t,r)) - \nabla F(\bar{X}(t))\bigr)$, $H_2(t)=\mathbf{B} \mathbf{Q}^{\top} \sum_{r=0}^{\tau-1} \bigl(\nabla F(\Phi(t,r)) - \nabla F(\bar{X}(t))\bigr) - \tau \mathbf{Q}^{\top} \bigl(\nabla F(\bar{X}(t+1)) - \nabla F(\bar{X}(t))\bigr)$, $\mathbf{B}=\mathbf{I}-\mathbf{\Lambda}$, $\mathbf{P}$ is a permutation matrix, and $\iota = \sqrt{-1}$ denotes the imaginary unit.
\end{lemma}
\begin{proof}
See Appendix \ref{lemma_most_important_proof}.
\end{proof}
The vector $D(t)$ characterizes the consensus deviations of $X(t)$ and $E(t)$ in the transformed coordinate system (\ref{transformed_system}). To see this, note that for any vector $X(t)$, the matrix $\mathbf{Q}$ satisfies $\mathbf{Q}\mathbf{Q}^{\top} = \mathbf{I} - \frac{1}{N}\mathbf{1}\mathbf{1}^{\top}$ and $\mathbf{Q}^{\top}\mathbf{Q} = \mathbf{I}$ from Lemma~\ref{lemma_from_LED1}. Consequently, we have
\begin{equation}\label{averaging_terms}
\left\{
\begin{aligned}
\|\mathbf{Q}^{\top} X(t)\|^2 &= \|\mathbf{Q}\mathbf{Q}^{\top} X(t)\|^2 = \|X(t) - \bar{X}(t)\|^2,\\
\|\mathbf{Q}^{\top} E(t)\|^2 &= \|\mathbf{Q}\mathbf{Q}^{\top} E(t)\|^2 = \|E(t) - \bar{E}(t)\|^2,
\end{aligned}
\right.
\end{equation}
where $\bar{E}(t) = \mathbf{1}_N \otimes \bar{e}(t)$. Since $D(t)$ is defined in terms of $\mathbf{Q}^{\top} X(t)$ and $\mathbf{Q}^{\top} E(t)$, its evolution governed by the contraction matrix $\Delta$ directly implies the exponential decay of the consensus errors of both $X(t)$ and $E(t)$. Thus, (\ref{transformed_system}) simplifies the convergence analysis of LPDVR to bounding $S(t)$ and $H(t)$, which originate from the gradient estimation error. The following section is devoted to deriving the required bounds via the properties of the gradient estimator.

\section{Preliminaries for Convergence Analysis}\label{Preliminaries_for_CA}

This section establishes auxiliary results for the convergence analysis. We first analyze the statistical properties of the gradient estimator in \eqref{gradient_process_algorithm1}, including its unbiasedness and variance bounds. We then investigate the evolution of the transformed consensus deviation $D(t)$, whose recursion in \eqref{transformed_system} is key to controlling the consensus error. These results will be instrumental in proving the main convergence theorems in the next sections.

\subsection{Expectation of Gradient Estimation}

In LPDVR, the randomness of $g_i(\cdot)$ arises from the independent random variables $\{ s_i(t,r)\}^{t\geq 0, 0\leq r <\tau}_{1\leq i\leq N}$. Let $\mathcal{F}_{t,r-1}$ denote the history of the system up to the $(r-1)$-th local update at outer iteration $t$, i.e., the $\sigma$-algebra generated by $\{s_i(t,k)\}_{1\leq i\leq N,\, 0\leq k<r}$ over all previous outer iterations. Owing to the independence of $\{s_i(t,r)\}$, given $\mathcal{F}_{t,r-1}$, the gradient estimator $g_i(\phi_i(t,r))$ is unbiased:
\begin{align}\label{expectation_gradient}
\mathbb{E}\left[g_i(\phi_i(t,r)) \mid \mathcal{F}_{t,r-1}\right] = \nabla f_i(\phi_i(t,r)),
\end{align}
for any $i\in\mathcal{V}$, $t\geq 0$, and $0\leq r<\tau$. Having established the unbiasedness of the gradient estimator, we next proceed to analyze its variance, which plays a critical role in bounding the error terms $S(t)$ and $H(t)$ in \eqref{transformed_system}.

\subsection{Variance of Gradient Estimation}
In this subsection, we analyze the variance of the gradient estimator $g_i(\phi_i(t,r))$. To this end, we introduce the following variables:
\begin{equation}\label{variance_definition}
\left\{
\begin{aligned}
\delta_i(t,r)=&\frac{1}{m_i}\sum^{m_i}_{h=1}\Vert z_{i,h}(t,r)-\bar{x}(t)\Vert^2,\\
\delta(t,r)=&\sum^{N}_{i=1}\delta_i(t,r),\quad \delta(t)=\sum^{\tau-1}_{r=0}\delta(t,r).
\end{aligned}
\right.
\end{equation}
For any $i\in\mathcal{V}$, $t\geq 0$, and $0\leq r<\tau$, the gradient estimation $g_i(\phi_i(t,r))$ generated in Algorithm \ref{algorithm_recur} satisfies
\begin{align}
&\mathbb{E}\big[\Vert g_i(\phi_i(t,r))-\nabla f_i(\phi_i(t,r))\Vert^2|\mathcal{F}_{t,r-1}\big]\nonumber\\
=&\mathbb{E}\big[\Vert \nabla f_{i,s_i(t,r)}(\phi_i(t,r))-\nabla f_{i,s_i(t,r)}(z_{i,s_i(t,r)}(t,r))\nonumber\\
       &-(\nabla f_i(\phi_i(t,r))-\frac{1}{m_i}\sum^{m_i}_{h=1}\nabla f_{i,h}(z_{i,h}(t,r)))\Vert^2|\mathcal{F}_{t,r-1}\big]\nonumber\\
\leq & \frac{1}{m_i}\sum^{m_i}_{h=1}\Vert \nabla f_{i,h}(\phi_i(t,r))-\nabla f_{i,h}(z_{i,h}(t,r))\Vert^2,\label{variance_step1}
\end{align}
where the equality follows from \eqref{gradient_process_algorithm1}, and the inequality follows from the standard conditional variance decomposition and the uniform sampling of $s_i(t,r)$. Furthermore, we have
\begin{align}
&\mathbb{E}\big[\Vert g_i(\phi_i(t,r))-\nabla f_i(\phi_i(t,r))\Vert^2\big]\nonumber\\
\leq & 2L^2\mathbb{E}\big[\Vert \phi_i(t,r)-\bar{x}(t)\Vert^2+\frac{1}{m_i}\sum^{m_i}_{h=1}\Vert \bar{x}(t)-z_{i,h}(t,r)\Vert^2\big]\nonumber\\
=&2L^2\mathbb{E}\big[\Vert \phi_i(t,r)-\bar{x}(t)\Vert^2\big]+2L^2 \mathbb{E}\big[\delta_i(t,r)\big],\label{variance_step2}
\end{align}
where the inequality and equality follow from  \eqref{variance_step1} and Assumption~\ref{smooth_assumption}, and  \eqref{variance_definition}, respectively. From the definitions (\ref{matrix_definition}), (\ref{matrix_definition2}), and (\ref{variance_definition}), summing \eqref{variance_step2} over $i$ and $r$ yields
\begin{align}
\mathbb{E}\big[ \sum^{\tau-1}_{r=0}\big\Vert S(t,r)\Vert^2\big]\leq 2L^2\mathbb{E}\big[\Vert \widehat{\Phi}(t)\Vert^2\big]+2L^2 \mathbb{E}\big[\delta(t)\big],\label{variance_step3}
\end{align}
where we define
\begin{align}\label{matrix_definition_average}
   \Vert \widehat{\Phi}(t)\Vert^2\triangleq\sum^{\tau-1}_{r=0} \Vert \Phi(t,r)-\bar{X}(t)\Vert^2
\end{align}
for notational convenience. We further define 
\begin{align}\label{definition_of_Mm}
M=\max_{1\leq i\leq N}\{m_i\},\quad\quad m=\min_{1\leq i\leq N}\{m_i\}.
\end{align} 
As shown in \eqref{variance_step3}, the variance of the gradient estimator $g_i(\phi_i(t,r))$ depends on $\Vert \widehat{\Phi}(t)\Vert^2$ and $\delta(t)$, whose properties are established below.
\begin{lemma}\label{variance_bound_importance1}
If the stepsize satisfies 
\begin{align}\label{stepsize_condition1}
\alpha\leq\min\Big\{\frac{1}{2\tau L},\frac{1}{5\sqrt{M}L}, \frac{1}{30L}\sqrt{\frac{m}{\tau M}}\Big\},
\end{align}
we have
\begin{equation}\label{end_part1}
\left\{
\begin{aligned}
&\mathbb{E}\big[\delta(t)\big]\leq M\mathbb{E}\big[\Vert X(t)-\bar{X}(t)\Vert^2\big]+\frac{M}{m}\mathbb{E}\big[\Vert \widehat{\Phi}(t)\Vert^2 \big],\\
&\mathbb{E}[\Vert \widehat{\Phi}(t)\Vert^2]\leq 40\alpha^2\tau^3 N\mathbb{E}[\Vert{\nabla f}(\bar{x}(t))\Vert^2]\\
&\quad\quad\quad\quad\quad\quad+212\tau\mathbb{E}[\Vert {D}(t)\Vert^2].
\end{aligned}
\right.
\end{equation}
\end{lemma}
\begin{proof}
See Appendix \ref{variance_bound_importance1_proof}.
\end{proof}
With Lemma~\ref{variance_bound_importance1} in place, it remains to bound $\mathbb{E}[\|D(t)\|^2]$. Recall that $D(t)$ governs both the consensus dynamics via the matrix $\Delta$ (see Lemma \ref{lemma_from_LED1}) and the variance of the gradient estimator as shown in Lemma~\ref{variance_bound_importance1}. Therefore, analyzing the evolution of $D(t)$ is important to analyze the convergence of LPDVR, which will be presented in the next subsection.

\subsection{Contractive Inequality for $D(t)$}

We now analyze the evolution of the transformed consensus deviation $D(t)$, which plays a central role in the convergence analysis. Combining Lemma~\ref{lemma_from_LED1} with the unbiasedness of the gradient estimator in \eqref{expectation_gradient} yields
\begin{align}
\mathbb{E}\big[\Vert {D}(t+1)\Vert^2\big]=&\mathbb{E}\big[\Vert \Delta{D}(t)-\alpha H(t)\Vert^2\big]+\alpha^2\mathbb{E}\big[\Vert S(t)\Vert^2\big]\nonumber\\
&+2\alpha^2\mathbb{E}\big[\langle S(t), H(t)\rangle\big].\label{second_part_2}
\end{align}
Using (\ref{second_part_2}) and $\Vert a+b\Vert^2\leq \frac{1}{\delta}\Vert a\Vert^2+\frac{1}{1-\delta}\Vert b\Vert^2$ gives
\begin{align}
\mathbb{E}\big[\Vert {D}(t+1)\Vert^2\big]\leq&\delta \mathbb{E}\big[\Vert {D}(t)\Vert^2\big]+\frac{2\alpha^2}{1-\delta}\mathbb{E}\big[\Vert {H}(t)\Vert^2\big]\nonumber\\
&+2\alpha^2\mathbb{E}\big[\Vert S(t)\Vert^2\big].\label{second_part_3}
\end{align}
Thus, to analyze the property of $\mathbb{E}[\Vert {D}(t+1)\Vert^2]$, it remains to bound  $\mathbb{E}[\Vert S(t)\Vert^2]$ and $\mathbb{E}[\Vert H(t)\Vert^2]$ in (\ref{second_part_3}). 

We first bound $\mathbb{E}\big[\Vert S(t)\Vert^2\big]$. From (\ref{S_property}) in Lemma \ref{additional_lemma} in Appendix~\ref{app:supporting_lemmas}, we have
\begin{align}
\Vert S(t)\Vert^2\leq2\tau\sum^{\tau-1}_{r=0}\Vert S(t,r)\Vert^2,\label{second_part_4}
\end{align}
where $S(t,r)=G(\Phi(t,r))-\nabla F(\Phi(t,r))$. Then, we have
\begin{align}
\mathbb{E}\big[\Vert S(t,r)&\Vert^2\big]=\sum^{N}_{i=1}\mathbb{E}\big[\Vert g_i(\phi_i(t,r))-\nabla f_i(\phi_i(t,r))\Vert^2\big]\nonumber\\
\leq&2L^2\mathbb{E}\big[\Vert \Phi(t,r)-\bar{X}(t)\Vert^2\big]+2L^2\mathbb{E}\big[\delta(t,r)\big],\label{second_part_5}
\end{align}
where the equality and the inequality follow from (\ref{matrix_definition}) and (\ref{variance_step2}), respectively. Furthermore, we have
\begin{align}
&\mathbb{E}\big[\Vert S(t)\Vert^2\big]\leq 4\tau L^2\sum^{\tau-1}_{r=0}\mathbb{E}\big[\Vert \Phi(t,r)-\bar{X}(t)\Vert^2+\delta(t,r)\big]\nonumber\\
&\quad\leq \frac{8\tau L^2 M}{m}\big(\mathbb{E}\big[\Vert\widehat{\Phi}(t)\Vert^2\big]+\frac{m}{2}\mathbb{E}\big[\Vert X(t)-\bar{X}(t)\Vert^2\big]\big),\label{second_part_6}
\end{align}
where the first inequality is obtained by substituting \eqref{second_part_4} into \eqref{second_part_5}, and the second follows from \eqref{delta_step8}.

Next, we bound $\mathbb{E}[\Vert H(t)\Vert^2]$. From the definition $\bar{s}(t,r)=\frac{1}{N}\sum_{i=1}^{N}\{g_i(t,r)-\nabla f_i(\phi_i(t,r))\}$, we have
\begin{align}
\Vert \sum^{\tau-1}_{r=0}\bar{s}(t,r)\Vert^2\leq \frac{\tau}{N}\sum^{\tau-1}_{r=0}\Vert S(t,r)\Vert^2.\label{second_part_7}
\end{align}
Furthermore, combining (\ref{second_part_4}), (\ref{second_part_6}), and (\ref{second_part_7}) yields
\begin{align}
\mathbb{E}\big[\Vert \sum^{\tau-1}_{r=0}\bar{s}(t,r)\Vert^2\big]\leq& \frac{2\tau L^2M}{N}\mathbb{E}\big[\Vert X(t)-\bar{X}(t)\Vert^2\big]\nonumber\\
&+\frac{4\tau L^2 M}{mN}\mathbb{E}\big[\Vert\widehat{\Phi}(t)\Vert^2\big].\label{second_part_8}
\end{align}
For $\mathbb{E}[\Vert H(t)\Vert^2]$ in (\ref{second_part_3}), its property is summarized in Lemma~\ref{additional_lemma}. Substituting (\ref{second_part_8}) into (\ref{second_part_9}) in Lemma \ref{additional_lemma} yields
\begin{align}
\mathbb{E}\big[\Vert H(t&)\Vert^2\big]\leq\big(4\tau L^2+\frac{8\alpha^2\tau^3L^4M}{m}\Vert \mathbf{B}^{-1}\Vert\big)\mathbb{E}\big[\Vert\widehat{\Phi}(t)\Vert^2\big]\nonumber\\
&+2\alpha^2\tau^2L^2\Vert \mathbf{B}^{-1}\Vert N\mathbb{E}\big[\Vert\sum^{\tau-1}_{r=0}\overline{\nabla F}(\Phi(t,r))\Vert^2\big]\nonumber\\
&+8\alpha^2\tau^3L^4\Vert \mathbf{B}^{-1}\Vert  { M}\mathbb{E}\big[\Vert X(t)-\bar{X}(t)\Vert^2\big].\label{second_part_10}
\end{align}
Finally, substituting $\Vert\mathbf{B}^{-1}\Vert=(1-\lambda)^{-1}$, (\ref{D_X_property}) in Lemma \ref{additional_lemma}, (\ref{second_part_6}), and (\ref{second_part_10}) into (\ref{second_part_3}) yields
\begin{align}
\mathbb{E}\big[&\Vert {D}(t+1)\Vert^2\big]\leq\frac{4\alpha^4\tau^2L^2 N}{(1-\delta)(1-\lambda)}\mathbb{E}\big[\Vert\sum^{\tau-1}_{r=0}\overline{\nabla F}(\Phi(t,r))\Vert^2\big]\nonumber\\
&+\Big(\delta +32\tau\alpha^2 L^2M+\frac{64\alpha^4\tau^3L^4  { M}}{(1-\delta)(1-\lambda)}\Big)\mathbb{E}\big[\Vert {D}(t)\Vert^2\big]\nonumber\\
&+\Big(\frac{24\tau\alpha^2 L^2 M}{m(1-\delta)}+\frac{16\alpha^4\tau^3L^4M}{m(1-\delta)(1-\lambda)}\Big)\mathbb{E}\big[\Vert\widehat{\Phi}(t)\Vert^2\big].\label{second_part_11}
\end{align}
Upon substitution of the second inequality in (\ref{end_part1}) into (\ref{second_part_11}) and simplification by collecting like terms, we arrive at
\begin{align}
&\mathbb{E}\big[\Vert {D}(t+1)\Vert^2\big]\leq\frac{4\alpha^4\tau^2L^2 N}{(1-\delta)(1-\lambda)}\mathbb{E}\big[\Vert\sum^{\tau-1}_{r=0}\overline{\nabla F}(\Phi(t,r))\Vert^2\big]\nonumber\\
&+\frac{80MN(12\tau^4\alpha^4 L^2 +8\alpha^6\tau^6L^4\Vert \mathbf{B}^{-1}\Vert) }{m(1-\delta)} \mathbb{E}\big[\Vert{\nabla f}(\bar{x}(t))\Vert^2\big]\nonumber\\
&+\big( 32\tau\alpha^2 L^2M+\frac{3456\alpha^4\tau^4L^4  { M}}{(1-\delta)(1-\lambda)}+\frac{5088\tau^2\alpha^2 L^2 M}{m(1-\delta)}\nonumber\\
&\quad\quad\ +\delta\big)\mathbb{E}\big[\Vert {D}(t)\Vert^2\big].\label{second_part_12}
\end{align}
If a stepsize satisfies
\begin{align}\label{stepsize_rule1}
\alpha\leq\sqrt{\frac{1-\lambda}{m}}\frac{1}{\tau L},
\end{align}
then (\ref{second_part_12}) can be simplified to
\begin{align}
&\mathbb{E}\big[\Vert {D}(t+1)\Vert^2\big]\leq \frac{1920\tau^4\alpha^4 L^2 MN}{m(1-\delta)} \mathbb{E}\big[\Vert{\nabla f}(\bar{x}(t))\Vert^2\big]\nonumber\\
&\quad+\Big(\delta +32\tau\alpha^2 L^2M+\frac{10176\tau^2\alpha^2 L^2 M}{m(1-\delta)}\Big)\mathbb{E}\big[\Vert {D}(t)\Vert^2\big] \nonumber\\
&\quad+\frac{4\alpha^4\tau^2L^2 N}{(1-\delta)(1-\lambda)}\mathbb{E}\big[\Vert\sum^{\tau-1}_{r=0}\overline{\nabla F}(\Phi(t,r))\Vert^2\big].\label{property_DT}
\end{align}

In summary, under the stepsize conditions (\ref{stepsize_condition1}) and (\ref{stepsize_rule1}), the transformed consensus deviation $D(t)$ satisfies the contraction-type inequality \eqref{property_DT}. This bound characterizes the evolution of the consensus error in terms of $\|\nabla f(\bar{x}(t))\|^2$ and  $\|\sum_{r=0}^{\tau-1}\overline{\nabla F}(\Phi(t,r))\|^2$. Moreover, with a sufficiently small stepsize, the contraction coefficient can be made less than unity, i.e., $\delta +32\tau\alpha^2 L^2M+\frac{10176\tau^2\alpha^2 L^2 M}{m(1-\delta)} < 1.$ This inequality serves as the basis for establishing the convergence rates of LPDVR for nonconvex, convex, and strongly convex cases in the subsequent sections.

\section{Convergence Analysis for Nonconvex Cases}\label{convergence_for_nonconvex_case}
Building on the contraction inequality and the gradient estimator properties established in Section \ref{Preliminaries_for_CA}, we now establish the convergence rates of LPDVR for the nonconvex case. To measure convergence, we introduce the following metric:
\begin{align}
\mathcal{E}(t) \triangleq \mathbb{E}[\|\nabla f(\bar{x}(t))\|^2] + \frac{1}{\tau}\sum_{r=0}^{\tau-1}\mathbb{E}[\|\frac{1}{N}\sum_{i=1}^{N} \nabla f_i(\phi_i(t,r))\|^2],\label{nonconvex_step6}
\end{align}
where the average $\bar{x}(t)=\frac{1}{N}\sum^{N}_{i=1}x_i(t)$ was defined in \eqref{averaging_x}. This quantity captures both the stationarity of the averaged optimization variable and the local gradient deviations across agents and inner iterations, which together characterize the convergence behavior of LPDVR. The following theorem gives the nonconvex convergence result for LPDVR.
\begin{theorem}\label{thm_nonconvex}
Under Assumptions \ref{smooth_assumption}, \ref{non_empty}, \ref{mixing_matrix}, and stepsize $\alpha\leq\mathcal{O}(\frac{1-\lambda}{\tau L})$ (see \eqref{stepsize_condition1}, \eqref{stepsize_rule1}, \eqref{stepsize_condition3}, \eqref{stepsize_condition4}, \eqref{stepsize_condition5}), it holds that
\begin{align*}
\frac{1}{T} \sum_{t=0}^{T-1} \mathcal{E}(t)
\leq \mathcal{O}\Big( \frac{f(\bar{x}(0)) - f(x^*)}{\alpha \tau T}\Big) + \mathcal{O}\Big(\frac{\alpha^2 \tau^2 L^2 \zeta_0^2 M}{(1-\lambda)^2 T} \Big),
\end{align*}
where $\zeta_0^2 \triangleq \frac{1}{N} \sum_{i=1}^N \|\nabla f_i(\bar{x}(0)) - \nabla f(\bar{x}(0))\|^2.$
\end{theorem}

\begin{proof}
The proof consists of three main steps. First, we derive a descent inequality for $\mathbb{E}[f(\bar{x}(t))]$ with a consensus error term.
Second, we use the recursion for $D(t)$ to bound the cumulative
consensus error in terms of $\|D(0)\|^2$ and
$\sum_{t=0}^{T-1}\mathcal{E}(t)$.
Finally, we sum the descent inequality and use the bound from Step 2
to obtain the claimed convergence rate under the stepsize conditions.

\textbf{Step 1: Descent inequality for $f(\bar{x}(t))$.} 
From the matrix forms (\ref{matrix_form_algorithm_x}), (\ref{matrix_form_algorithm_y}) and the definitions in (\ref{matrix_definition2}), we can obtain the following property
\begin{align}
\hspace{-0.18cm}\bar{x}(t+1)=\bar{x}(t)-\frac{\alpha}{N}\sum^{\tau-1}_{r=0}\sum^{N}_{i=1}\big( \nabla f_i(\phi_i(t,r))+s_i(t,r)\big),\label{nonconvex_step1}
\end{align}
where $s_i(t,r)=g_i(t,r)-\nabla f_i(\phi_i(t,r))$. From Assumption \ref{smooth_assumption}, Lemma \ref{smooth_property_nes},  (\ref{expectation_gradient}), and (\ref{nonconvex_step1}), we can directly obtain  the inequality 

\noindent $\mathbb{E}\big[f(\bar{x}(t+1))\big] \leq {\alpha^2 \tau L}\mathbb{E}\big[\sum^{\tau-1}_{r=0}\Vert \frac{1}{N}\sum^{N}_{i=1} \nabla f_i(\phi_i(t,r))\Vert^2\big]+\mathbb{E}\big[f(\bar{x}(t))\big]-\frac{\alpha}{N}\mathbb{E}\big[\big\langle \nabla f(\bar{x}(t)),\sum^{\tau-1}_{r=0}\sum^{N}_{i=1} \nabla f_i(\phi_i(t,r))\big\rangle\big]+ {\alpha^2 L}\mathbb{E}\big[\Vert \sum^{\tau-1}_{r=0}\bar{s}(t,r)\Vert^2\big].$ Furthermore, substituting (\ref{second_part_8}) into this inequality yields
\begin{align}
&\mathbb{E}\big[f(\bar{x}(t+1))\big]\leq {\alpha^2 \tau L}\mathbb{E}\big[\sum^{\tau-1}_{r=0}\big\Vert \frac{1}{N}\sum^{N}_{i=1} \nabla f_i(\phi_i(t,r))\big\Vert^2\big]\nonumber\\
&+\mathbb{E}\big[f(\bar{x}(t))\big]-{\alpha}\mathbb{E}\big[\big\langle \nabla f(\bar{x}(t)),\frac{1}{N}\sum^{\tau-1}_{r=0}\sum^{N}_{i=1} \nabla f_i(\phi_i(t,r))\big\rangle\big]\nonumber\\
&+\frac{2\alpha^2\tau L^3M}{N}\big(\frac{2}{m}\mathbb{E}\big[\Vert \widehat{\Phi}(t)\Vert^2\big]+\mathbb{E}\big[\Vert X(t)-\bar{X}(t)\Vert^2\big]\big).\label{nonconvex_step2}
\end{align}
For $-{\alpha}\big\langle \nabla f(\bar{x}(t)),\frac{1}{N}\sum^{\tau-1}_{r=0}\sum^{N}_{i=1} \nabla f_i(\phi_i(t,r))\big\rangle$ in (\ref{nonconvex_step2}), we have
\begin{align}
&-{\alpha}\big\langle \nabla f(\bar{x}(t)),\frac{1}{N}\sum^{\tau-1}_{r=0}\sum^{N}_{i=1} \nabla f_i(\phi_i(t,r))\big\rangle\leq \frac{\alpha L^2}{2N}\Vert \widehat{\Phi}(t)\Vert^2\nonumber\\
&-\frac{\alpha}{2}\sum^{\tau-1}_{r=0}\Vert \frac{1}{N}\sum^{N}_{i=1} \nabla f_i(\phi_i(t,r))\Vert^2 -\frac{\alpha\tau}{2}\Vert \nabla f(\bar{x}(t))\Vert^2,\label{nonconvex_step3}
\end{align}
which follows from Assumption \ref{smooth_assumption}, (\ref{matrix_definition}), and (\ref{matrix_definition_average}). Substituting (\ref{end_part1}), (\ref{nonconvex_step3}), and (\ref{D_X_property}) in Lemma \ref{additional_lemma} into (\ref{nonconvex_step2}) gives
\begin{align}
\mathbb{E}&\big[f(\bar{x}(t+1))\big]\leq  \mathbb{E}\big[f(\bar{x}(t))\big]-\frac{\alpha\tau}{2}\mathbb{E}\big[\Vert \nabla f(\bar{x}(t))\Vert^2\big]\nonumber\\
&+\Big(\frac{106\tau\alpha L^2}{N}+\frac{848\tau^2\alpha^2 L^3 M}{mN}+\frac{8\alpha^2\tau L^3M}{N}\Big)\mathbb{E}\big[\Vert {D}(t)\Vert^2\big]\nonumber\\
&+\Big({20\alpha^3\tau^3 L^2}+\frac{160\alpha^4\tau^4 L^3 M}{m}\Big) \mathbb{E}\big[\Vert{\nabla f}(\bar{x}(t))\Vert^2\big]\nonumber\\
&-\frac{\alpha}{2}\big(1-2\alpha L\tau\big)\sum^{\tau-1}_{r=0}\mathbb{E}\big[\Vert \frac{1}{N}\sum^{N}_{i=1} \nabla f_i(\phi_i(t,r))\Vert^2\big].\label{nonconvex_step4}
\end{align}
If a stepsize satisfies
\begin{align}\label{stepsize_condition3}
\alpha\leq\min\Big\{\frac{1}{14\tau L},\frac{m}{9M\tau L}\Big\},
\end{align}
we have $\frac{160\alpha^4\tau^4 L^3 M}{m}\leq{20\alpha^3\tau^3 L^2}$ and $\frac{848\tau^2\alpha^2 L^3 M}{mN}\leq\frac{106\tau\alpha L^2}{N}$. Furthermore, (\ref{nonconvex_step4}) can be simplified to 
\begin{align}
\mathbb{E}\big[f(\bar{x}(t&+1))\big]\leq  \mathbb{E}\big[f(\bar{x}(t))\big]-\frac{\alpha\tau}{4}\mathbb{E}\big[\Vert \nabla f(\bar{x}(t))\Vert^2\big]\nonumber\\
&-\frac{\alpha}{4}\sum^{\tau-1}_{r=0}\mathbb{E}\big[\Vert \frac{1}{N}\sum^{N}_{i=1} \nabla f_i(\phi_i(t,r))\Vert^2\big]\nonumber\\
&+\big(\frac{212\tau\alpha L^2}{N}+\frac{8\alpha^2\tau L^3M}{N}\big)\mathbb{E}\big[\Vert {D}(t)\Vert^2\big].\label{add_more_for_similarity_low3}
\end{align}
After rearranging terms of (\ref{add_more_for_similarity_low3}), summing over $t=0,1,\cdots,T-1$, and dividing by $\frac{\alpha\tau}{4T}$, we obtain
\begin{align}
\frac{1}{T}\sum^{T-1}_{t=0}\mathcal{E}(t)\leq  &\frac{848 L^2+32\alpha L^3M}{NT}\sum^{T-1}_{t=0}\mathbb{E}[\Vert {D}(t)\Vert^2]\nonumber\\
&+\frac{4(f(\bar{x}(0))-f(x^*))}{\alpha\tau T},\label{nonconvex_step7}
\end{align}
where $\mathcal{E}(t)$ is defined in (\ref{nonconvex_step6}).

\textbf{Step 2: Analysis of the consensus deviation $D(t)$.} 
Combining (\ref{property_DT}) and (\ref{nonconvex_step6}) gives
\begin{align}
&\mathbb{E}\big[\Vert {D}(t+1)\Vert^2\big]\leq\frac{1920\tau^4\alpha^4 L^2 MN}{m(1-\sqrt{\lambda})(1-\lambda)} \mathcal{E}(t)\nonumber\\
&+\Big(\delta +32\tau\alpha^2 L^2M+\frac{10176\tau^2\alpha^2 L^2 M}{m(1-\sqrt{\lambda})}\Big)\mathbb{E}\big[\Vert {D}(t)\Vert^2\big].\label{nonconvex_step8}
\end{align}
Because of $\lambda<1$, if we design the stepsize satisfies
\begin{align}\label{stepsize_condition4}
    \alpha\leq\min\left\{\frac{(1-\lambda)\sqrt{m}}{500\tau L\sqrt{M}},\frac{(1-\sqrt{\lambda})\sqrt{m}}{16\sqrt{3}\tau L\sqrt{M}},\frac{\sqrt{1-\sqrt{\lambda}}}{8L\sqrt{2\tau M}}\right\},
\end{align}
we can obtain the following inequalities
\begin{equation}\label{inequalities_add}
\left\{
\begin{aligned}
\delta +32\tau\alpha^2 L^2M+\frac{10176\tau^2\alpha^2 L^2 M}{m(1-\delta)}&\leq \frac{1+\delta}{2}\triangleq\bar{\delta},\\
\frac{1920\tau^4\alpha^4 L^2 M}{m(1-\delta)(1-\lambda)}&\leq \tau^2\alpha^2.
\end{aligned}
\right.
\end{equation}
Moreover, substituting (\ref{inequalities_add}) into (\ref{nonconvex_step8}) gives $\mathbb{E}\big[\Vert {D}(t+1)\Vert^2\big]\leq\bar{\delta}\mathbb{E}\big[\Vert {D}(t)\Vert^2\big]+{\tau^2\alpha^2N} \mathcal{E}(t),$ which further implies
\begin{align}\label{add_more_for_similarity_low1}
\mathbb{E}\big[\Vert {D}(t)\Vert^2\big]\leq\bar{\delta}^{t}\Vert {D}(0)\Vert^2+{\tau^2\alpha^2N} \sum^{t-1}_{l=0}\bar{\delta}^{t-1-l}\mathcal{E}(l).
\end{align}
Summing \eqref{add_more_for_similarity_low1} over $0\leq t\leq T-1$ and dividing by $T$ yields $\frac{1}{T}\sum^{T-1}_{t=0}\mathbb{E}[\Vert {D}(t)\Vert^2]\leq\frac{\Vert {D}(0)\Vert^2}{T}\sum^{T-1}_{t=0}\bar{\delta}^{t}+\frac{\tau^2\alpha^2N}{T}\sum^{T-1}_{t=0} \sum^{t-1}_{l=0}\bar{\delta}^{t-1-l}\mathcal{E}(l).$ Thus, applying Lemma~\ref{ming_yan_lemma} simplifies this inequality to
\begin{align}
\frac{1}{T}\sum^{T-1}_{t=0}\mathbb{E}\big[\Vert {D}(t)\Vert^2\big]
\leq\frac{\Vert {D}(0)\Vert^2}{(1-\bar{\delta})T}+\frac{\tau^2\alpha^2N}{(1-\bar{\delta})T}\sum^{T-1}_{t=0} \mathcal{E}(t).\label{nonconvex_step9}
\end{align}
Substituting (\ref{nonconvex_step9}) into (\ref{nonconvex_step7}) and rearranging terms gives
\begin{align}
\hspace{-0.65em}&\Big(1-\frac{848 L^2\tau^2\alpha^2+32  L^3\tau^2\alpha^3 M}{1-\bar{\delta}}\Big)\frac{1}{T}\sum^{T-1}_{t=0}\mathcal{E}(t)\nonumber\\
\hspace{-0.65em}\leq&  \frac{4(f(\bar{x}(0))-f(x^*))}{\alpha\tau T}+ \frac{\big({848 L^2}+{32\alpha L^3M}\big)\Vert {D}(0)\Vert^2}{(1-\bar{\delta})NT}.\label{nonconvex_step10}
\end{align}

\textbf{Step 3: Combining the bounds.} 
With $\bar{\delta}=\frac{1+\delta}{2}$, if the stepsize satisfies
\begin{align}\label{stepsize_condition5}
\alpha\leq\min\Big\{\frac{1-\sqrt{\lambda}}{30\sqrt{10}\tau L},\frac{\sqrt[3]{1-\sqrt{\lambda}}}{7\tau L\sqrt[3]{M}}\Big\},
\end{align}
we have $\frac{848 L^2\tau^2\alpha^2}{1-\bar{\delta}}\leq \frac{1}{4}$ and $\frac{32 L^3\tau^2\alpha^3 M}{1-\bar{\delta}}\leq \frac{1}{4}$. Thus,  (\ref{nonconvex_step10}) can be simplified to $\frac{1}{T}\sum^{T-1}_{t=0}\mathcal{E}(t)\leq   \frac{({3392 L^2}+{128\alpha L^3M})\Vert {D}(0)\Vert^2}{(1-\sqrt{\lambda})NT}+\frac{8(f(\bar{x}(0))-f(x^*))}{\alpha\tau T}.$ Because of $x_i(0)=\mathbf{0}_n$ and $y_i(0)=\mathbf{0}_n$, we have $\|{D}(0)\|^2 
\leq \|{\mathbf{B}}^{-\frac{1}{2}} \mathbf{Q}^{\top} \mathbf{E}(0)\|^2 \leq \frac{2 \alpha^2 \tau^2}{1 - \lambda} \|\nabla F(\bar{{X}}(0)) - \mathbf{1} \otimes \nabla f({\bar{x}}(0))\|^2.$ Thus, the proof is completed.
\end{proof}
Theorem~\ref{thm_nonconvex} establishes a sublinear convergence rate for LPDVR in the nonconvex setting. The following observations can be made:
\begin{itemize}
    \item First, the metric $\mathcal{E}(t)$ defined in \eqref{nonconvex_step6} captures both the gradient norm at the averaged iterate and the average of local gradients across agents and inner iterations. Moreover, under the convergence result in Theorem~\ref{thm_nonconvex}, the consensus deviation $D(t)$ satisfies the bound in \eqref{nonconvex_step9}, implying that consensus is achieved in the time-averaged sense at the rate $\mathcal{O}(1/T)$.
    \item Second, the $\mathcal{O}\big(\frac{f(\bar{x}(0)) - f(x^*)}{\alpha \tau T}\big)$ term reflects the standard descent behavior, and the $\mathcal{O}\big(\frac{\alpha^2 \tau^2 L^2 \zeta_0^2{M}}{(1-\lambda)^2 T}\big)$ term depends on the network connectivity through $(1-\lambda)^{-2}$ and the number of local samples through ${M}$. These dependencies have intuitive interpretations. Networks with stronger connectivity (where $\lambda$ is closer to $0$) facilitate more rapid information diffusion, thereby accelerating convergence. Conversely, a larger local sample size $M$ introduces greater uncertainty into the gradient estimation established via the variance reduction technique (\ref{gradient_process_algorithm1}), which naturally weakens the resulting convergence bound.
\end{itemize}

In addition to the convergence result in Theorem~\ref{thm_nonconvex},
the following corollary provides a bound on the time-averaged
consensus error of LPDVR in the nonconvex setting.
\begin{corollary}[Consensus in the Nonconvex Setting]\label{corollary_consensus_nonconvex}
Under the same assumptions as Theorem~\ref{thm_nonconvex}, consensus of LPDVR is achieved in the time-averaged sense at a rate of $\mathcal{O}(1/T)$, i.e., $\frac{1}{T} \sum_{t=0}^{T-1}\sum_{i=1}^{N} \mathbb{E}[\|x_i(t)-\bar{x}(t)\|^2] \leq \mathcal{O}({1}/{T}).$
\end{corollary}

Corollary \ref{corollary_consensus_nonconvex} follows directly from Theorem~\ref{thm_nonconvex}, the bound on $\mathbb{E}[\|D(t)\|^2]$ established in \eqref{nonconvex_step9}, together with the relation $\|X(t)-\bar{X}(t)\|^2 \leq 4\|D(t)\|^2$ from Lemma~\ref{additional_lemma}. Therefore, although Theorem~\ref{thm_nonconvex} primarily characterizes convergence of the averaged iterate $\bar{x}(t)$ to a stationary point, the above result further ensures that all agents reach consensus in the time-averaged sense at the same rate. In other words, the algorithm not only drives the average trajectory to a stationary point, but also guarantees that the individual agents' optimization variables converge to a common stationary point.

Having established the convergence result for nonconvex objectives, we now turn to the convex and strongly convex cases, where the analysis follows the same three-step framework but measures progress via the distance to the optimal solution.

\section{Convergence Analysis of Convex and Strongly Convex Cases}\label{convergence_analysis_for_c_and_sc}
We now consider the case where each local function $f_i$ is convex or strongly convex. The following assumption is imposed throughout this section.

\begin{assumption}[Convexity]\label{assumption_convexity}
For each $i \in \mathcal{V}$, the local function $f_i: \mathbb{R}^n \to \mathbb{R}$ is $\mu$-strongly convex with $0 \leq \mu \leq L$. In the case where $\mu = 0$, each $f_i$ reduces to a standard convex function.
\end{assumption}

We begin by deriving a unified inequality that will serve as the foundation for both convex and strongly convex analyses.

\begin{lemma}\label{lemma_convex_unified}
Under Assumptions \ref{smooth_assumption}, \ref{non_empty}, \ref{mixing_matrix} and \ref{assumption_convexity}, and that the stepsize satisfies \eqref{stepsize_condition1} and \eqref{convex_stepsize1}, then we have
\begin{align}
&\hspace{-0.19cm}\mathbb{E}[\Vert \bar{x}(t+1)-x^*\Vert^2]\leq\frac{1484\alpha L\tau+16\tau\alpha^2 L^2M}{N}\mathbb{E}[\Vert {D}(t)\Vert^2]\nonumber\\
&\hspace{-0.19cm}+(1-\frac{\alpha\tau\mu}{2})\mathbb{E}[\Vert \bar{x}(t)-x^*\Vert^2]-\frac{\alpha\tau}{2}\mathbb{E}[f(\bar{x}(t))-f(x^*)].\hspace{-0.2cm}\label{convex_step8}
\end{align}
\end{lemma}

\begin{proof}
From (\ref{expectation_gradient}) and (\ref{nonconvex_step1}), we have
\begin{align}
\mathbb{E}\big[\Vert \bar{x}(t&+1)-x^*\Vert^2\big]\leq2\alpha^2\mathbb{E}\big[\Vert \frac{1}{N}\sum^{\tau-1}_{r=0}\sum^{N}_{i=1} \nabla f_i(\phi_i(t,r))\Vert^2\big]\nonumber\\
&+\mathbb{E}\big[\Vert \bar{x}(t)-x^*\Vert^2\big]+2\alpha^2\mathbb{E}\big[\Vert \sum^{\tau-1}_{r=0}\bar{s}(t,r)\Vert^2\big]\nonumber\\
&-\frac{2\alpha}{N}\mathbb{E}\big[\big\langle \bar{x}(t)-x^*,\sum^{\tau-1}_{r=0}\sum^{N}_{i=1} \nabla f_i(\phi_i(t,r))\big\rangle\big].\label{convex_step1}
\end{align}
For the term $-\frac{2\alpha}{N}\mathbb{E}[\langle \bar{x}(t)-x^*,\sum^{\tau-1}_{r=0}\sum^{N}_{i=1} \nabla f_i(\phi_i(t,r))\rangle]$ in (\ref{convex_step1}), we have
\begin{align}
&-\frac{2\alpha}{N}\sum^{\tau-1}_{r=0}\sum^{N}_{i=1}\big\langle \bar{x}(t)-x^*, \nabla f_i(\phi_i(t,r))\big\rangle\leq \frac{2\alpha L}{N}\big\Vert \widehat{\Phi}(t)\big\Vert^2\nonumber\\
&\quad-{2\alpha\tau}\big(f(\bar{x}(t))-f(x^*)+\frac{\mu}{4}\big\Vert \bar{x}(t)-x^*\big\Vert^2\big),\label{convex_step3}
\end{align}
which follows from Lemma \ref{scaffold_lemma} and (\ref{matrix_definition_average}). For the term $\Vert \frac{1}{N}\sum^{\tau-1}_{r=0}\sum^{N}_{i=1} \nabla f_i(\phi_i(t,r))\Vert^2$ in (\ref{convex_step1}), we have
\begin{align}
&\Vert \frac{1}{N}\sum^{\tau-1}_{r=0}\sum^{N}_{i=1} \nabla f_i(\phi_i(t,r))\Vert^2\nonumber\\
\leq &\frac{2\tau L^2}{N}\sum^{\tau-1}_{r=0}\sum^{N}_{i=1}\Vert  \phi_i(t,r)-\bar{x}(t)\Vert^2+2\tau^2\Vert \nabla f(\bar{x}(t))\Vert^2\nonumber\\
= &\frac{2\tau L^2}{N}\Vert  \widehat{\Phi}(t)\Vert^2+2\tau^2\Vert \nabla f(\bar{x}(t))\Vert^2,\label{convex_step5}
\end{align}
where the inequality and equality follow from Assumption \ref{smooth_assumption} and (\ref{matrix_definition2}), respectively. Substituting (\ref{second_part_8}),  (\ref{convex_step3}), and (\ref{convex_step5}) into (\ref{convex_step1}) and applying the condition $\alpha L\tau \leq 1$ together with (\ref{D_X_property}) in Lemma \ref{additional_lemma} yields the following result:
\begin{align}
&\mathbb{E}\big[\Vert \bar{x}(t+1)-x^*\Vert^2\big]\leq\big(1-\frac{\alpha\tau\mu}{2}\big)\mathbb{E}\big[\Vert \bar{x}(t)-x^*\Vert^2\big]\nonumber\\
&+4\alpha^2\tau^2\mathbb{E}[\Vert \nabla f(\bar{x}(t))\Vert^2]+\frac{8\tau \alpha^2L^2 M+6\alpha Lm}{mN}\mathbb{E}[\Vert \widehat{\Phi}(t)\Vert^2]\nonumber\\
&+\frac{16\tau\alpha^2 L^2M}{N}\mathbb{E}[\Vert {D}(t)\Vert^2]-{2\alpha\tau}\mathbb{E}\big[f(\bar{x}(t))-f(x^*)\big].\label{convex_step6}
\end{align}
If we choose the stepsize such that
\begin{align}
\alpha\leq\min\Big\{\frac{m}{8\tau LM},\frac{1}{380\tau L}\Big\}\label{convex_stepsize1}
\end{align}
which implies $ \frac{8\tau \alpha^2L^2 M}{mN}\leq\frac{\alpha L}{N}$ and $ {2\alpha\tau}-568\alpha^2\tau^2L\geq \frac{\alpha\tau}{2}$. Thus, together with (\ref{end_part1}) and the property $\Vert\nabla f(\bar{x}(t)) \Vert^2\leq 2L(f(\bar{x}(t))-f(x^*))$ from Assumption \ref{smooth_assumption}, then  (\ref{convex_step6}) can be simplified to (\ref{convex_step8}). The proof of Lemma \ref{lemma_convex_unified} is completed.
\end{proof}

\subsection{Convex Case ($\mu = 0$)}
For the convex case, we measure the convergence in terms of the optimality gap:
\begin{align}\label{definition_convex_E}
\mathcal{E}(t) \triangleq \mathbb{E}[f(\bar{x}(t))] - f(x^*).
\end{align}
The following theorem gives the convex convergence result for LPDVR.
\begin{theorem}\label{thm_convex}
Under Assumptions \ref{smooth_assumption}, \ref{non_empty}, \ref{mixing_matrix}, \ref{assumption_convexity} ($\mu=0$), and stepsize $\alpha\leq\mathcal{O}(\frac{1-\lambda}{\tau L})$ (see \eqref{stepsize_condition1},  \eqref{stepsize_rule1}, \eqref{convex_stepsize1},  \eqref{stepsize_condition_lemma}, and \eqref{convex_stepsize3}), it holds that
\begin{align*}
\frac{1}{T} \sum_{t=0}^{T-1} \mathcal{E}(t)
\leq \mathcal{O}\Big( \frac{\Vert \bar{x}(0) - x^*\Vert^2}{\alpha \tau T} \Big)+ \mathcal{O}\Big(\frac{\alpha^2 \tau^2 L \zeta_0^2{M}}{(1-\lambda)^2 T} \Big),
\end{align*}
where $\zeta_0^2 \triangleq \frac{1}{N} \sum_{i=1}^N \|\nabla f_i(\bar{x}(0)) - \nabla f(\bar{x}(0))\|^2.$
\end{theorem}

\begin{proof}
The proof, which follows the same strategy as that of
Theorem~\ref{thm_nonconvex}, consists of three main steps.

\textbf{Step 1: Descent inequality for $\|\bar{x}(t)-x^*\|^2$.} 
 Substituting $\mu=0$ into (\ref{convex_step8}) and rearranging terms yields 
 \begin{align}
\mathcal{E}(t)\leq&\frac{2}{\alpha\tau}\big(\mathbb{E}[\Vert \bar{x}(t)-x^*\Vert^2]-\mathbb{E}[\Vert \bar{x}(t+1)-x^*\Vert^2]\big)\nonumber\\
&+\big(\frac{2968 L}{N}+\frac{32\alpha L^2M}{N}\big)\mathbb{E}[\Vert {D}(t)\Vert^2].\label{add_for_tac_pagelimitation1}
 \end{align}
Then, summing (\ref{add_for_tac_pagelimitation1}) over $t=0,\ldots,T-1$ gives
\begin{align}
&\frac{1}{T}\sum^{T-1}_{t=0}\mathcal{E}(t)\leq\frac{2\Vert \bar{x}(0)-x^*\Vert^2}{\alpha\tau T}\nonumber\\
&\quad\quad\quad+\frac{2968 L+32\alpha L^2M}{N T}\sum^{T-1}_{t=0}\mathbb{E}\big[\Vert {D}(t)\Vert^2\big].\label{convex_step10}
\end{align}

\textbf{Step 2: Analysis of the consensus deviation $D(t)$.} 
By substituting (\ref{convex_step5}) into (\ref{property_DT}), we can obtain  
\begin{align}
\mathbb{E}[\Vert {D}(t&+1)\Vert^2]\leq\frac{8\alpha^4\tau^3L^4 }{(1-\delta)(1-\lambda)}\big(\mathbb{E}[\Vert  \widehat{\Phi}(t)\Vert^2]\nonumber\\
&+\frac{241\tau   MN}{mL^2} \mathbb{E}[\Vert{\nabla f}(\bar{x}(t))\Vert^2]\big)+\big(\delta+32\tau\alpha^2 L^2M\nonumber\\
&+\frac{10176\tau^2\alpha^2 L^2 M}{m(1-\delta)} \big)\mathbb{E}[\Vert {D}(t)\Vert^2].\label{add_property_for_convex}
\end{align}
If a stepsize satisfies
\begin{align}\label{stepsize_condition_lemma}
\alpha\leq\min\Big\{\frac{(1-\lambda)\sqrt{m}}{300\tau L\sqrt{M}},\frac{\sqrt{1-\sqrt{\lambda}}}{8\sqrt{2\tau M}L}\Big\},
\end{align}
together with $\Vert \nabla f(\bar{x}(t))\Vert^2\leq 2L(f(\bar{x}(t))-f(x^*))$ from Assumption \ref{smooth_assumption}, (\ref{end_part1}) and (\ref{definition_convex_E}), (\ref{add_property_for_convex}) can be simplified to
\begin{align}
\mathbb{E}\big[\Vert {D}(t+1)\Vert^2\big]\leq\tilde{\delta}\mathbb{E}\big[\Vert {D}(t)\Vert^2\big]+4\alpha^2\tau^2 LN\mathcal{E}(t),\label{convex_step11}
\end{align}
where $\tilde{\delta}=\frac{3+\delta}{4}$. Let $d_t=\mathbb{E}[\|D(t)\|^2]$.
Summing \eqref{convex_step11} over $t=0,\ldots,T-1$
and rearranging terms yields  $(1-\tilde{\delta})\sum_{t=0}^{T-1}d_t
\leq d_0-d_T
+4\alpha^2\tau^2LN\sum_{t=0}^{T-1}\mathcal{E}(t).$
Since $d_T\geq0$ and
$1-\tilde{\delta}=(1-\delta)/4$,
dividing both sides by $(1-\tilde{\delta})T$
gives
\begin{align}
\hspace{-0.25cm}\frac{1}{T}\sum^{T-1}_{t=0}\mathbb{E}[\Vert {D}(t)\Vert^2]\leq \frac{4\Vert {D}(0)\Vert^2}{(1-{\delta})T} +\frac{16\alpha^2\tau^2 LN}{(1-\delta)T}\sum^{T-1}_{t=0}\mathcal{E}(t).\hspace{-0.11cm}\label{add_property_for_convex2}
\end{align}

\textbf{Step 3: Combining the bounds.} 
Substituting (\ref{add_property_for_convex2}) into (\ref{convex_step10}) yields
\begin{align}
&\big(1-\frac{47488\alpha^2\tau^2 L^2}{1-\delta}-\frac{512\alpha^3\tau^2 L^3 M}{1-\delta}\big)\frac{1}{T}\sum^{T-1}_{t=0}\mathcal{E}(t)\nonumber\\
\leq& \frac{2\Vert \bar{x}(0)-x^*\Vert^2}{\alpha\tau T}+\frac{4(2968 L+32\alpha L^2M)\Vert {D}(0)\Vert^2}{N (1-{\delta})T}.\label{convex_step13}
\end{align}
If a stepsize satisfies
\begin{align}\label{convex_stepsize3}
\alpha\leq\min\Big\{\frac{\sqrt{1-\sqrt{\lambda}}}{380\tau L},\frac{\sqrt[3]{1-\sqrt{\lambda}}}{13L\sqrt[3]{M\tau^2}}\Big\},
\end{align}
we have $ \frac{47488\alpha^2\tau^2 L^2}{1-\delta}\leq\frac{1}{3}$ and $\frac{512\alpha^3\tau^2 L^3 M}{1-\delta}\leq\frac{1}{3}$. Moreover, because $x_i(0)=\mathbf{0}_n$ and $y_i(0)=\mathbf{0}_n$, we have $\|{D}(0)\|^2 
\leq  \frac{2 \alpha^2 \tau^2}{1 - \lambda} \|\nabla F(\bar{{X}}(0)) - \mathbf{1} \otimes \nabla f({\bar{x}}(0))\|^2.$ Substituting these inequalities into (\ref{convex_step13}) completes the proof.
\end{proof}

Theorem~\ref{thm_convex} establishes a sublinear convergence rate of $\mathcal{O}(1/T)$ for LPDVR in the convex setting. The first term in the bound captures the standard descent behavior of the averaged iterate, while the second term reflects the influence of the initial consensus error and the network connectivity through $(1-\lambda)^{-2}$. The interpretation of these dependencies is similar to the nonconvex case: better-connected networks yield faster convergence, while a larger number of local samples introduces more gradient estimation uncertainty. The key distinction lies in the convergence metric: the convex analysis measures the optimality gap \eqref{definition_convex_E}, whereas the nonconvex analysis tracks the gradient norm and local gradient deviations \eqref{nonconvex_step6}. This difference reflects the distinct notions of stationarity in the two settings. Notably, the stepsize condition $\alpha = \mathcal{O}((1-\lambda)/(\tau L))$ remains unchanged, confirming that the variance reduction mechanism does not require a more conservative stepsize when moving from nonconvex to convex settings. 

In addition to the convergence result in Theorem~\ref{thm_nonconvex},
the following corollary provides a bound on the time-averaged
consensus error of LPDVR in the convex setting.

\begin{corollary}[Consensus in the Convex Setting]\label{corollary_consensus_convex}
Under the same assumptions as Theorem~\ref{thm_convex}, consensus of LPDVR is achieved in the time-averaged sense at a rate of $\mathcal{O}(1/T)$, i.e., $\frac{1}{T} \sum_{t=0}^{T-1}\sum_{i=1}^{N} \mathbb{E}[\|x_i(t)-\bar{x}(t)\|^2] \leq \mathcal{O}({1}/{T}).$
\end{corollary}
Corollary \ref{corollary_consensus_convex} follows directly from Theorem~\ref{thm_convex}, the bound on $\mathbb{E}[\|D(t)\|^2]$ in \eqref{add_property_for_convex2}, and the relation $\|X(t)-\bar{X}(t)\|^2 \leq 4\|D(t)\|^2$ in Lemma~\ref{additional_lemma}. Similar to the nonconvex case, the convex result not only ensures convergence of the averaged iterate to the optimal solution, but also guarantees that all agents reach consensus in the time-averaged sense at the same rate. Building on the sublinear results for both nonconvex and convex settings, we next turn to the strongly convex case.

\subsection{Strongly Convex Case ($\mu > 0$)}
For the strongly convex case, the following theorem establishes a linear convergence rate for LPDVR.
\begin{theorem}\label{thm_strongly_convex}
Suppose the stepsize $\alpha\leq\mathcal{O}(\frac{\mu(1-\lambda)}{\tau L^2})$ (see \eqref{stepsize_condition1},  \eqref{stepsize_rule1}, \eqref{convex_stepsize1},  \eqref{stepsize_condition_lemma}, and (\ref{stepsize_condition_strongly_convex})). Then the iterates generated by Algorithm~\ref{algorithm_recur} satisfy
\begin{align}
\mathbb{E}[\|\bar{x}(t) - x^*\|^2] \leq \mathcal{O}\Big(\big(1-\frac{\alpha\tau\mu}{4}\big)^{t}\Big).
\end{align}
\end{theorem}

\begin{proof}
The proof consists of three main steps.
First, we derive a descent inequality for
$\mathbb{E}[\|\bar{x}(t)-x^*\|^2]$.
Second, we obtain a recursion for the consensus error.
Finally, we combine the two inequalities and bound the spectral
radius of the resulting system to establish linear convergence.

\textbf{Step 1: Descent inequality for $\|\bar{x}(t)-x^*\|^2$.} 
From $f(x^*)\leq f(\bar{x}(t))$, (\ref{convex_step8}) can be simplified to
\begin{align}
\mathbb{E}\big[\Vert \bar{x}(t&+1)-x^*\Vert^2\big]\leq\big(1-\frac{\alpha\tau\mu}{2}\big)\mathbb{E}\big[\Vert \bar{x}(t)-x^*\Vert^2\big]\nonumber\\
&\quad+\frac{1484\alpha L\tau+16\tau\alpha^2 L^2M}{N}\mathbb{E}\big[\Vert {D}(t)\Vert^2\big].\label{strongly_convex_step1}
\end{align}

\textbf{Step 2: Analysis of the consensus deviation $D(t)$.} If a stepsize satisfies (\ref{stepsize_condition_lemma}), substituting $\Vert \nabla f(\bar{x}(t))\Vert^2\leq L^2\Vert \bar{x}(t)-x^*\Vert^2$ from Assumption \ref{smooth_assumption} into (\ref{convex_step11}) and dividing both sides by $N$ yields 
\begin{align}
&\frac{1}{N}\mathbb{E}\big[\Vert {D}(t+1)\Vert^2\big]\nonumber\\
\leq&\frac{3+\delta}{4N}\mathbb{E}\big[\Vert {D}(t)\Vert^2\big]+2\alpha^2\tau^2L^2\mathbb{E}\big[\Vert\bar{x}(t)-x^*\Vert^2\big].\label{add2_forTAC}
\end{align}

\textbf{Step 3: Spectral radius of the coupled linear system.} 
Combining (\ref{strongly_convex_step1}) and (\ref{add2_forTAC}) yields a linear system
\begin{align}\label{system_matrix_form}
\begin{bmatrix}
\mathbb{E}\big[ \|\bar{x}(t+1) - x^*\|^2 \big] \\[6pt]
\frac{1}{N} \mathbb{E}\big[ \|{D}(t+1)\|^2 \big]
\end{bmatrix} \preceq A
\begin{bmatrix}
\mathbb{E}\big[ \|\bar{x}(t) - x^*\|^2 \big] \\[6pt]
\frac{1}{N} \mathbb{E}\big[ \|{D}(t)\|^2 \big]
\end{bmatrix},
\end{align}
where 
\begin{align}
    A=\begin{bmatrix}
1 - \dfrac{\alpha \tau \mu}{2} & 1484\alpha L\tau + 16\tau\alpha^2 L^2 M \\[8pt]
2\alpha^2 \tau^2 L^2 & \dfrac{3+\delta}{4}
\end{bmatrix}.
\end{align}
Thus, based on the definition of $\Vert\cdot\Vert_1$, we have $\rho(A)\leq\Vert A\Vert_1=\max\{1- \dfrac{\alpha \tau \mu}{2}+2\alpha^2 \tau^2 L^2,1484\alpha L\tau + 16\tau\alpha^2 L^2 M+   \dfrac{3+\delta}{4}\}.$ Thus, if the stepsize satisfies
\begin{align}\label{stepsize_condition_strongly_convex}
\alpha\leq\min\Big\{\frac{\mu}{8\tau L^2},\frac{1-\sqrt{\lambda}}{6400\tau L},\frac{\sqrt{1-\sqrt{\lambda}}}{32 L\sqrt{\tau M}}\Big\},
\end{align}
we have $\rho(A)<1-\frac{\alpha\tau\mu}{4}$. Taking the norm on both sides of \eqref{system_matrix_form} proves Theorem~\ref{thm_strongly_convex}.
\end{proof}

Theorem~\ref{thm_strongly_convex} establishes a linear convergence rate for LPDVR in the strongly convex setting. The rate is governed by the spectral radius $\rho(A)$ of the coupled system matrix in \eqref{system_matrix_form}, which captures the joint contraction of the optimality gap and the consensus deviation. Under the stepsize condition $\alpha = \mathcal{O}(\mu(1-\lambda)/(\tau L^2))$, the spectral radius satisfies $\rho(A) < 1 - \alpha\tau\mu/4$, yielding the advertised linear rate. This result is the first of its kind for variance-reduced decentralized optimization with multiple local updates under stochastic gradients.

\begin{corollary}[Consensus in the Strongly Convex Setting]\label{corollary_consensus_strongly_convex}
Under the same assumptions as Theorem~\ref{thm_strongly_convex}, consensus of LPDVR is achieved at a linear rate, i.e., $\mathbb{E}\big[\|x_i(t)-\bar{x}(t)\|^2\big] \leq\mathcal{O}\big((1-\frac{\alpha\tau\mu}{4})^{t}\big)$ for any $i\in\mathcal{V}$.
\end{corollary}
Corollary \ref{corollary_consensus_strongly_convex} follows from the fact that $\mathbb{E}[\|D(t)\|^2]$ and $\mathbb{E}[\|\bar{x}(t) - x^*\|^2]$ are governed by the same coupled linear system in \eqref{system_matrix_form}, whose spectral radius satisfies $\rho(A) < 1 - \alpha\tau\mu/4$ under the stepsize condition in Theorem~\ref{thm_strongly_convex}. Combined with the relation $\|X(t)-\bar{X}(t)\|^2 \leq 4\|D(t)\|^2$ in Lemma~\ref{additional_lemma}, we conclude that consensus is also achieved at a linear rate. Unlike the nonconvex and convex cases where consensus is guaranteed only in the time-averaged sense, the strongly convex setting ensures that the individual agent states converge to a common value at a linear rate, simultaneously with the convergence of the averaged iterate to the optimal solution.

\begin{remark}\label{remark_proof_structure}
The convergence analyses for the three cases (Theorems \ref{thm_nonconvex}, \ref{thm_convex}, and \ref{thm_strongly_convex}) share a unified three-step framework.
\begin{itemize}
\item In Step 1, we derive a recursive inequality. For the nonconvex case, this inequality relates $\mathbb{E}[f(\bar{x}(t+1))]$ to $\mathbb{E}[f(\bar{x}(t))]$, with additional terms involving the gradient norm, the local gradient deviations, and the consensus deviation. For the convex and strongly convex cases, an analogous inequality relates $\mathbb{E}[\|\bar{x}(t+1)-x^*\|^2]$ to $\mathbb{E}[\|\bar{x}(t)-x^*\|^2]$, with additional terms involving the optimality gap and the consensus deviation.
\item In Step 2, we characterize the contraction of the consensus deviation $D(t)$, whose recursion is governed by the contraction matrix $\Delta$. This contraction is expressed in terms of $\mathcal{E}(t)$ for the nonconvex and convex cases, while for the strongly convex case it is coupled with $\mathbb{E}[\|\bar{x}(t)-x^*\|^2]$.
\item In Step 3, we substitute the consensus deviation bound from Step 2 into the recursive inequality from Step 1, thereby eliminating the cumulative consensus error. This yields the final convergence results: a sublinear rate in terms of the gradient norm for the nonconvex case, a sublinear rate in terms of the optimality gap for the convex case, and a linear convergence rate for the strongly convex case.
\end{itemize}
\end{remark}

\begin{table*}[t]
\begin{center}
\begin{threeparttable}
\caption{Comparison of LPDVR and existing decentralized optimization algorithms that support multiple local updates. All methods are evaluated under \(\tau\) local updates per communication round, a network’s mixing rate \(\lambda = \Vert W - \frac{1}{N} \mathbf{1}_N \mathbf{1}_N^{\top} \Vert\), and a network of size $N$. }
\label{sample-table2}
\vskip -0.05in
\centering
\begin{small}
\setlength{\tabcolsep}{4pt}
\renewcommand{\arraystretch}{1.22}
\begin{tabular}{c c c c c c c }
\toprule
\multirow{2}{*}{\textbf{Reference}} & \textbf{Function}        & {\textbf{Convergence}}
& $\ $ \textbf{Stochastic} $\ $ & \textbf{Variance} $\ $ & \textbf{Communication} &  \textbf{Stepsize} \\
                                   & \textbf{Assumption}\tnote{1}      &             \textbf{Result}\tnote{2}           & \textbf{Gradient?}  & \textbf{Reduction?}  & \textbf{Overhead}\tnote{3}  & \textbf{Condition} \\
\midrule
\multirow{3}{*}{This Paper}
    &  $\mu$-s.c.  & Linear $\to 0$           & \multirow{3}{*}{\cmark} & \multirow{3}{*}{\cmark}  & \multirow{3}{*}{1} & $\mathcal{O}({\mu(1-\lambda)}/{\tau L^2})$ \\
    &  cvx      &  Sublinear $\to 0$         &                                  &                  &                          & $\mathcal{O}({(1-\lambda)}/{\tau L})$ \\
    &  noncvx     & Sublinear $\to 0$   &                                         &           &                          & $\mathcal{O}({(1-\lambda)}/{\tau L})$ \\
\midrule
\multirow{2}{*}{\cite{ren2025communication}\tnote{4}}      &  \multirow{2}{*}{noncvx}  & Sublinear $\propto \alpha\sigma^2$                 & \multirow{2}{*}{\cmark} &{\xmark} & \multirow{2}{*}{1}  & $\mathcal{O}(\lambda_l/\tau^2 L)$ \\
     &   & Sublinear $\to 0$                &  & \cmark &   & $\mathcal{O}(\lambda_l/\sqrt{N}\tau^3 L)$ \\
\midrule
\multirow{3}{*}{\cite{LED}}
    &  $\mu$-s.c.  & Linear $\propto \alpha\sigma^2$           & \multirow{3}{*}{\cmark} & \multirow{3}{*}{\xmark} & \multirow{3}{*}{1} & $\mathcal{O}({\mu(1-\lambda)}/{\tau L^2})$ \\
    &  cvx      & $\ $ Sublinear $\propto \alpha\sigma^2$         &    &                                                 &                          & $\mathcal{O}({(1-\lambda)}/{\tau L})$ \\
    &  noncvx     & Sublinear $\propto \alpha\sigma^2$    &            &                                         &                          & $\mathcal{O}({(1-\lambda)}/{\tau L})$ \\
\midrule
\cite{taolin1}              &  noncvx  & Sublinear $\propto \alpha\sigma^2$                 & \cmark & \xmark & 2  & $\mathcal{O}({(1-\lambda)}/{\tau L})$ \\
\cite{kunyuan1}  &  noncvx  & Sublinear $\to 0$                 & \xmark & \xmark & 2  & $\mathcal{O}({(1-\lambda)}/{\tau L})$ \\
\cite{LOCAL_SMALL_STEPSIZE} &  $\mu$-s.c.  & Linear $\to 0$        & \xmark & \xmark  & 2   &  $\mathcal{O}({\mu(1-\lambda)}/{\tau^2 L^2})$ \\
\cite{WeiShi3}                   & $\mu$-s.c.  & Linear $\to 0$                        & \xmark & \xmark & 2    & $\mathcal{O}({(1-\lambda)^2\sqrt\mu}/{\sqrt{NL}\tau^2 L})$ \\
\cite{wu2025effectiveness} & $\mu$-s.c.  & Linear $\to 0$                        & \xmark & \xmark  & 2    & $\mathcal{O}({\mu(1-\lambda)^2}/{\tau L^2})$ \\
\midrule
\multirow{2}{*}{\cite{liu2025guaranteeing}}      &  \multirow{2}{*}{noncvx}  & Sublinear $\propto \alpha\sigma^2$                & {\cmark} & \xmark & \multirow{2}{*}{1}  & \multirow{2}{*}{$\mathcal{O}({(1-\lambda)}/{\tau^4 L})$} \\
    &   & Sublinear $\to 0$                & {\xmark} & \xmark &   &  \\
   % \midrule
  %  \multirow{2}{*}{\cite{Mind_jieliu}}      &  \multirow{2}{*}{$\mu$-s.c.}  & Linear $\propto \alpha\sigma^2$                & {\cmark} & \xmark & \multirow{2}{*}{1}  & \multirow{2}{*}{$\mathcal{O}({\mu}/{\tau L^2})$} \\
  %  &   & Linear $\to 0$                & {\xmark} & \xmark &   &  \\
    \midrule
\multirow{2}{*}{\cite{huang2025beyond}}      & $\mu$-s.c.  & Linear $\propto \alpha\sigma^2$                 & \multirow{2}{*}{\cmark} & \multirow{2}{*}{\xmark} & \multirow{2}{*}{2}  &  $\mathcal{O}({\mu(1-\lambda)}/{\tau L^2})$ \\
     &  {noncvx}  & Sublinear $\propto \alpha\sigma^2$               &  & &  &  $\mathcal{O}({(1-\lambda)}/{\tau L})$  \\
     \midrule
\multirow{2}{*}{\cite{huang2025distributed}}      & noncvx+PL  & Linear $\propto \alpha\sigma^2$                 & \multirow{2}{*}{\cmark} & \multirow{2}{*}{\xmark} & \multirow{2}{*}{2}  &  $\mathcal{O}({\mu(1-\sqrt{1-\lambda^2})}/{\tau L^2})$ \\
     &  {noncvx}  & Sublinear $\propto \alpha\sigma^2$               &  & &  &  $\mathcal{O}({(1-\sqrt{1-\lambda^2})}/{\tau L})$  \\
\bottomrule
\end{tabular}
\end{small}
\vskip -0.05in
\begin{tablenotes}    
       \footnotesize  
       \item[{1}] {$\mu$-s.c., cvx, noncvx, and PL denote $\mu$-strongly convex, convex, nonconvex, and Polyak--Lojasiewicz condition, respectively. }
       \item[{2}] {Sublinear denotes $\mathcal{O}(1/T)$ convergence rate throughout this table, and \(\sigma^2\) denotes the variance of the stochastic gradient.}
        \item[{3}] {We measure communication overhead by the number of $n$-dimensional variables exchanged at communication rounds.}
        \item[{4}] {\cite{ren2025communication} is established over the graph Laplacian matrix, and $\lambda_l$ denotes its smallest nonzero eigenvalue.}
      \end{tablenotes} 
\end{threeparttable}
\end{center}
\vskip -0.15in
\end{table*}

\section{Comparisons with Existing Works}\label{comparisons_with_many_existing_works}
In this section, we compare LPDVR with representative existing decentralized optimization algorithms that support multiple local updates. Table~\ref{sample-table2} summarizes the comparison across five dimensions: function assumptions, convergence results, gradient type, communication overhead, and stepsize conditions. As is clearly shown in Table~\ref{sample-table2}, LPDVR offers several distinct advantages over existing methods, which are summarized below.

\begin{itemize}
    \item \textbf{Exact and Linear Convergence Guarantees.} By incorporating variance reduction, LPDVR achieves exact convergence in expectation under stochastic gradients, whereas most existing methods either converge only to a neighborhood \cite{LED,taolin1,liu2025guaranteeing,huang2025beyond,huang2025distributed} or rely on exact gradients to guarantee convergence \cite{kunyuan1,LOCAL_SMALL_STEPSIZE,WeiShi3,liu2025guaranteeing}. Moreover, LPDVR is the first decentralized algorithm with multiple local updates that achieves exact linear convergence for strongly convex objectives via variance reduction.

    \item \textbf{Improved Stepsize Design.} Although the ADMM-based variance-reduced method \cite{ren2025communication} achieves exact convergence in expectation for the nonconvex case, its stepsize scales as $\mathcal{O}(\lambda_l/(\sqrt{N}\tau^3 L))$, which becomes impractically small in large-scale networks due to the dependence on both $N$ and $\tau^3$. In contrast, the stepsize of LPDVR is independent of $N$ and depends only inversely on $\tau$, taking the form $\mathcal{O}((1-\lambda)/(\tau L))$ for both convex and nonconvex cases, which enables significantly larger stepsizes and faster convergence. Moreover, unlike the Laplacian-based design in \cite{ren2025communication}, our stepsize condition applies to any mixing matrix satisfying standard conditions, offering greater flexibility in network topology design. Notably, this stepsize condition is comparable to, or even less restrictive than, the stepsize conditions required by conventional stochastic or exact gradient methods \cite{liu2025guaranteeing,LOCAL_SMALL_STEPSIZE,kunyuan1,LED,taolin1,WeiShi3,huang2025beyond,huang2025distributed}, confirming that variance reduction is achieved without enforcing a more conservative stepsize.

   \item  \textbf{Communication Efficiency.} LPDVR requires each agent to transmit only its primal variable per communication round. In contrast, many existing decentralized algorithms require the exchange of two variables per round \cite{taolin1,LOCAL_SMALL_STEPSIZE,kunyuan1,WeiShi3,huang2025beyond,huang2025distributed}, typically including an auxiliary variable in addition to the optimization variable. Transmitting only one variable instead of two cuts the per-round communication overhead in half, making LPDVR particularly advantageous for bandwidth-limited environments and large-scale applications.
    
\end{itemize}

\begin{remark}
We emphasize that the communication strategy of LPDVR is fundamentally different from randomized communication-skipping approaches \cite{scaffnew,guo2026local,luyaoguo_skip_tcns,guanghongyang}. LPDVR performs a fixed number $\tau$ of local updates between rounds in a deterministic manner, whereas randomized methods skip communication randomly according to a certain probability. This distinction alters the algorithmic behavior and complicates its theoretical analysis. In LPDVR, the randomness arises solely from stochastic gradient estimates, not from communication skipping as in \cite{scaffnew,guo2026local,luyaoguo_skip_tcns,guanghongyang}. Therefore, the two lines of work address orthogonal sources of randomness and are not directly comparable. Rather, they complement each other in understanding different aspects of communication-efficient decentralized optimization.
\end{remark}

\section{Numerical Evaluations}\label{numerical_evaluations}

We consider the empirical risk minimization problem:
\begin{align}\label{simulation_euqation}
\hspace{-0.3cm}\min_{x\in\mathbb{R}^n}f(x)=\frac{1}{N}\sum^{N}_{i=1}\Big\{\frac{1}{m_i}\sum^{m_i}_{j=1}\Vert M_i x-b_{ij}\Vert^2+r_i\Vert x\Vert^2\Big\}.
\end{align}
This problem has broad applications, including distributed estimation \cite{wang_poor}, where each agent $i\in\mathcal{V}$ has $m_i$ noisy measurements $b_{ij} = M_i x + w_{ij}$ of an unknown parameter $x\in\mathbb{R}^n$, with $M_i\in\mathbb{R}^{n\times n}$ and $w_{ij}\in\mathbb{R}^{n}$ denoting the measurement matrix and noise, respectively. Each agent also has a regularization parameter, $r_i\ge 0$. The estimation task can be written as the empirical risk minimization in \eqref{simulation_euqation}. In our experiments, we take $M_i=\mathbf{I}_n$ and $r_i=1$ for all $i$. Thus, the local objective function reduces to
\(f_i(x)=\frac{1}{m_i}\sum^{m_i}_{j=1}(\Vert  x-b_{ij}\Vert^2+\Vert x\Vert^2),\)
which satisfies Assumptions~\ref{smooth_assumption} and~\ref{non_empty}.

\begin{figure}[htpb]
    \centering

    \begin{minipage}[t]{0.27\linewidth}
        \centering
        \includegraphics[width=\textwidth]{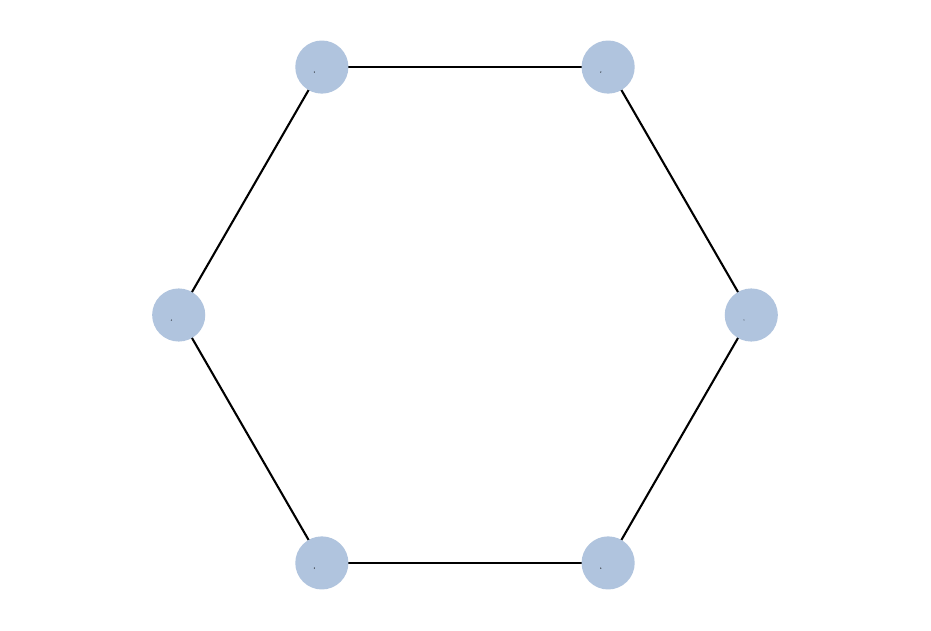}
        \\[2pt]
        {\small (a) Ring graph}
    \end{minipage}%
\hspace{2em}  
    \begin{minipage}[t]{0.27\linewidth}
        \centering
        \includegraphics[width=\textwidth]{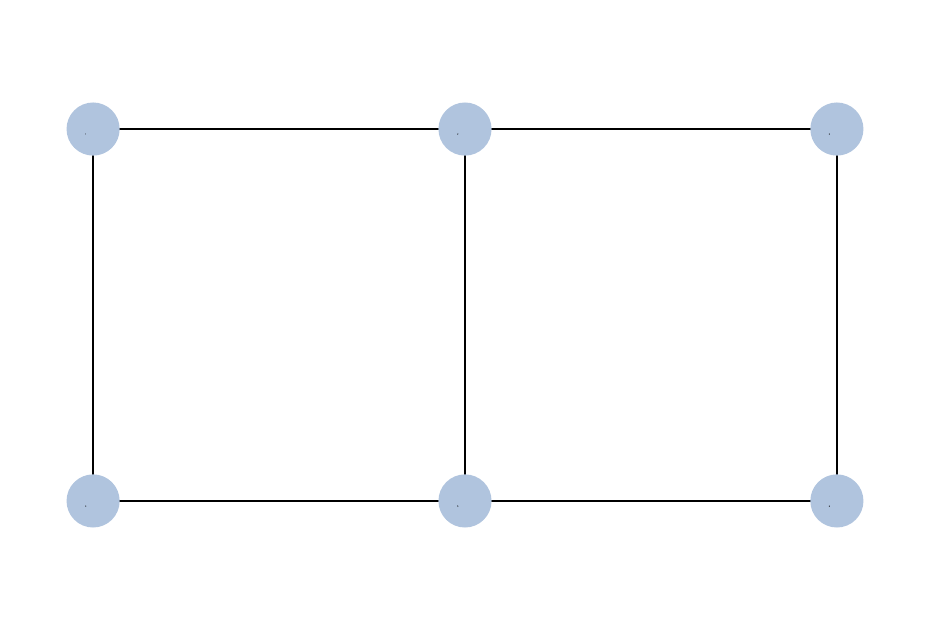}
        \\[2pt]
        {\small (b) Grid graph}
    \end{minipage}

    \caption{Ring and grid graphs with 6 nodes.}
    \label{fig:graph_topology}
     \vskip -0.3in
\end{figure}

We consider \(N\) agents, each with \(10\) noisy measurements (\(m_i=10\)), under two communication topologies: a ring graph and a \(2 \times (N/2)\) grid graph (see Fig.~\ref{fig:graph_topology}). For both topologies, the mixing matrix is constructed as \(W=\frac{1}{2}\mathbf{I}_N+\frac{1}{2}\widetilde{W}\), where \(\widetilde{W}\) is built according to the Metropolis rule \cite{xiao2006distributed}. Specifically, for each edge \((i,j)\), the weight is set as \([\widetilde{W}]_{ij} = 1/\max\{d_i, d_j\}\), with \(d_i\) denoting the degree of agent \(i\). Since \(\max\{d_i,d_j\}=\max\{d_j,d_i\}\), this rule yields a symmetric and doubly stochastic weight matrix. Each entry of \(b_{ij}\) is independently drawn from a uniform distribution over \([-10,10]\). We compare LPDVR with LT-ADMM-VR \cite{ren2025communication}, DIGing \cite{WeiShi3}, GTA \cite{LOCAL_SMALL_STEPSIZE}, and LED \cite{LED} on the above empirical risk minimization problem. The stepsizes of LT-ADMM-VR, GTA, and LED are chosen according to their respective theoretical guidelines (see Table \ref{sample-table2}). Since both DIGing \cite{WeiShi3} and GTA \cite{LOCAL_SMALL_STEPSIZE} are gradient-tracking-based algorithms, the stepsize of DIGing is set to be the same as that of GTA for a fair comparison. Convergence is measured by the distance from the averaged iterate to the optimal solution after each communication round: $e(k)=\Vert \bar{x}(k) - x^* \Vert$, where \(k\) denotes the communication round index. 
  \vskip -0.1in

\subsection{Performance Comparison with Varying Local Updates}

\begin{figure}[htpb]
    \begin{minipage}[t]{0.48\linewidth}
        \includegraphics[width=\textwidth]{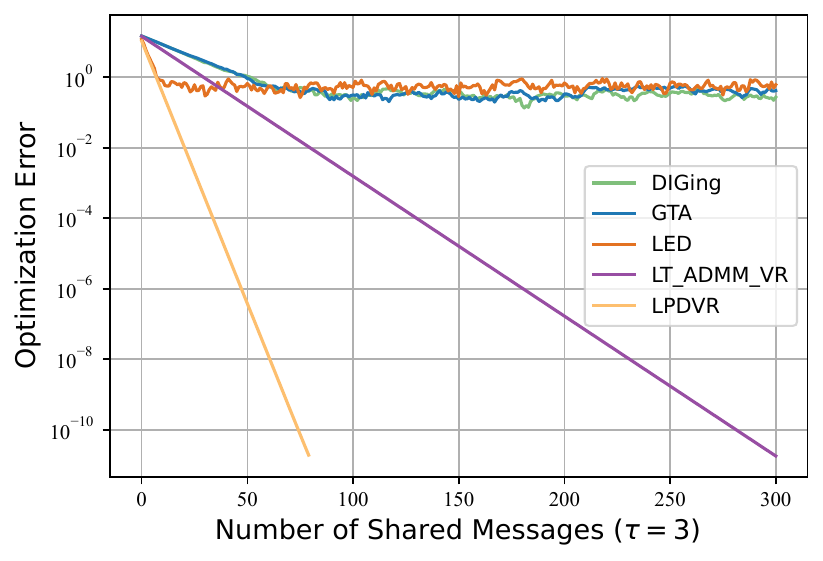}
    \end{minipage}   \hfill
    \begin{minipage}[t]{0.48\linewidth}
        \includegraphics[width=\textwidth]{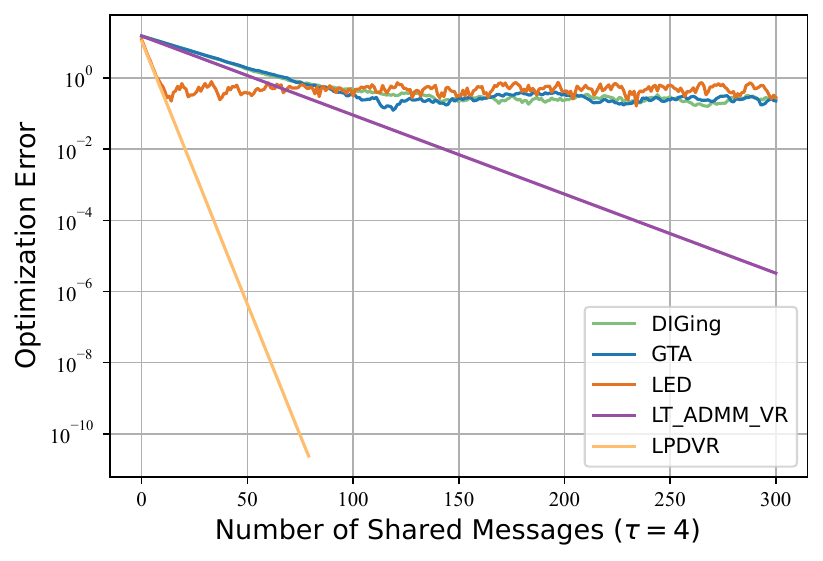}
    \end{minipage}   
    \begin{minipage}[t]{0.48\linewidth}
        \includegraphics[width=\textwidth]{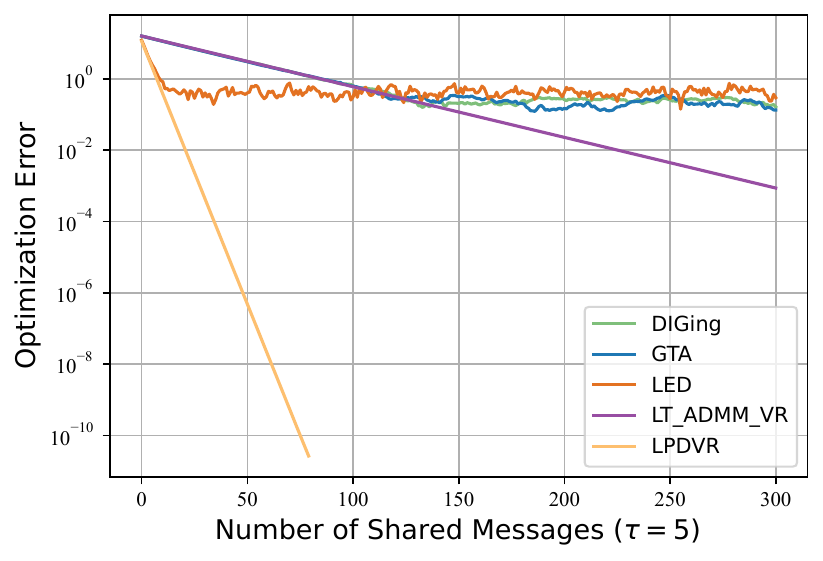}
    \end{minipage}\hfill
    \begin{minipage}[t]{0.48\linewidth}
        \includegraphics[width=\textwidth]{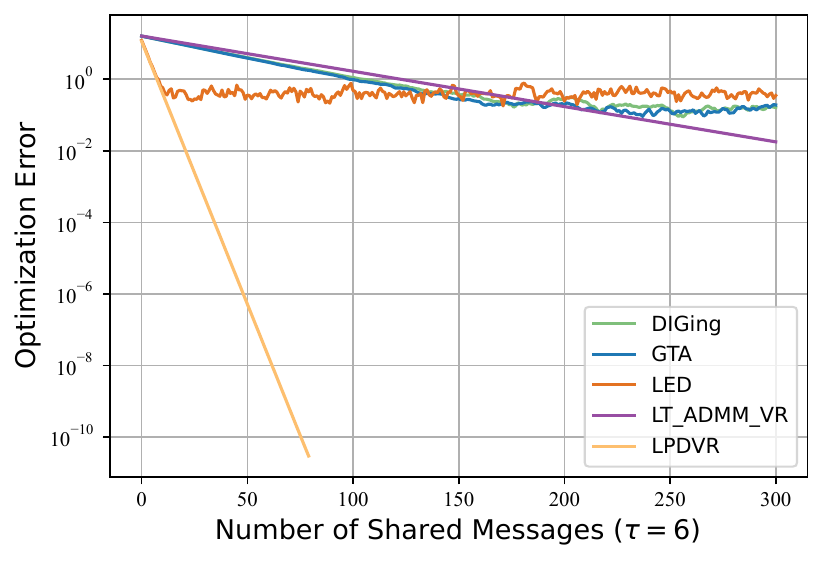}
    \end{minipage}
   \caption{Comparison of LPDVR with LT-ADMM-VR \cite{ren2025communication},
             DIGing \cite{WeiShi3}, GTA \cite{LOCAL_SMALL_STEPSIZE}, and LED \cite{LED} over the \textbf{ring} graph under different local updates $\tau=3,4,5,6$.}
    \label{ring_different_tau}
  \vskip -0.2in
\end{figure}

\begin{figure}[htpb]
    \begin{minipage}[t]{0.48\linewidth}
        \includegraphics[width=\textwidth]{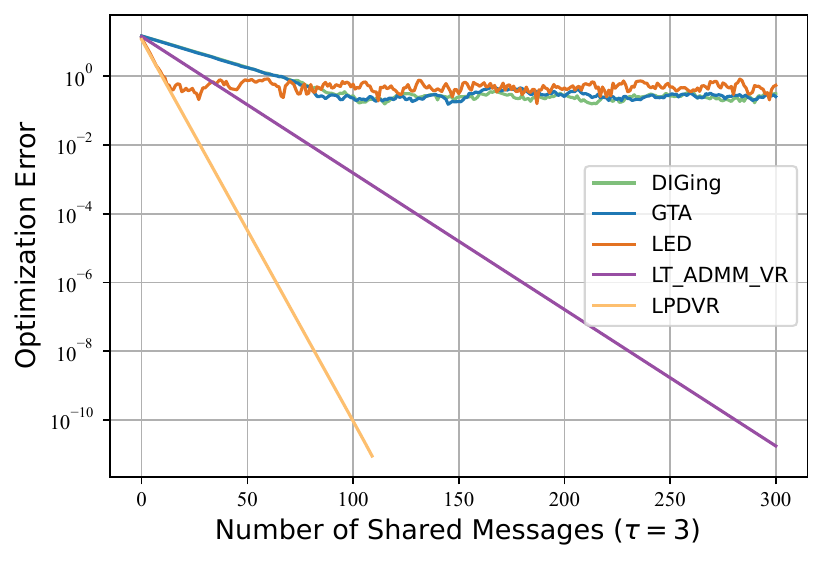}
    \end{minipage}   \hfill
    \begin{minipage}[t]{0.48\linewidth}
        \includegraphics[width=\textwidth]{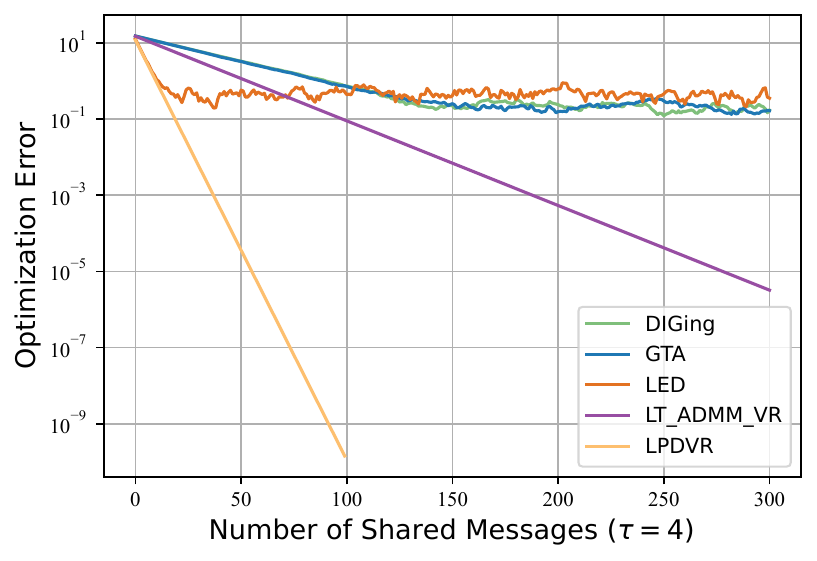}
    \end{minipage}    
    \begin{minipage}[t]{0.48\linewidth}
        \includegraphics[width=\textwidth]{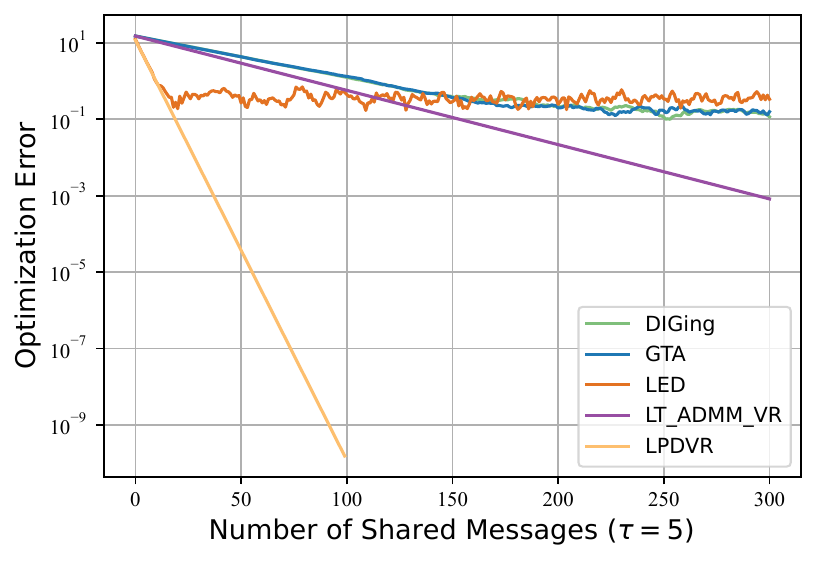}
    \end{minipage}\hfill
    \begin{minipage}[t]{0.48\linewidth}
        \includegraphics[width=\textwidth]{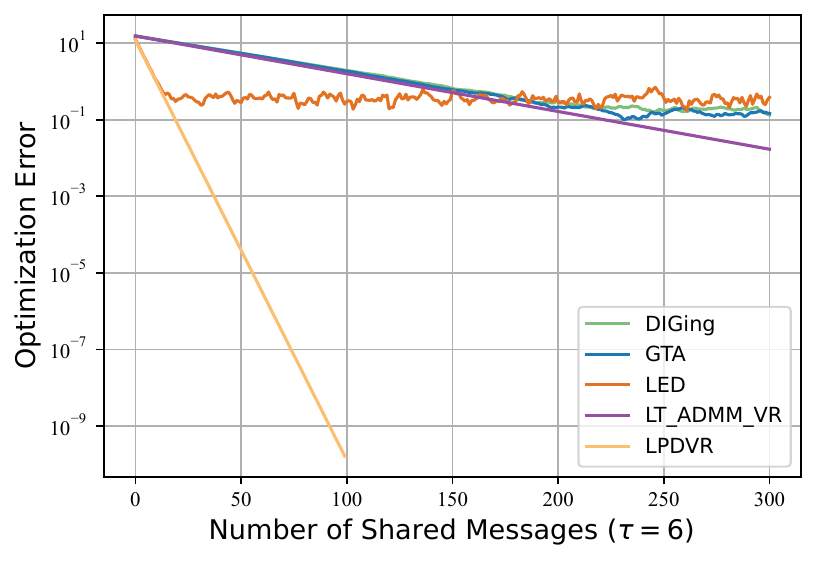}
    \end{minipage}
    \caption{Comparison of LPDVR with LT-ADMM-VR \cite{ren2025communication},
             DIGing \cite{WeiShi3}, GTA \cite{LOCAL_SMALL_STEPSIZE}, and LED \cite{LED}
              over  the \textbf{grid} graph under different local updates $\tau=3,4,5,6$.}
    \label{grid_different_tau}
    \vskip -0.2in
\end{figure}

We first investigate the convergence performance of LPDVR and competing algorithms under a fixed network size of \(N=6\) with varying numbers of local updates, \(\tau=3,4,5,6\). This experiment is conducted under both ring and grid topologies, and the two settings lead to consistent observations. As shown in Fig.~\ref{ring_different_tau}, under the ring topology, LPDVR and LT-ADMM-VR achieve significantly higher accuracy than DIGing, GTA, and LED, whose errors stagnate at a neighborhood of around \(10^{-1}\). Among the two high-accuracy methods, LPDVR converges significantly faster than LT-ADMM-VR across all tested values of \(\tau\). Similar observations can be made from Fig.~\ref{grid_different_tau} under the grid topology. Moreover, as \(\tau\) increases, LPDVR exhibits increasingly faster convergence, and the performance gap between LPDVR and LT-ADMM-VR becomes more pronounced under both topologies. This observation is corroborated by our theoretical analysis: LPDVR admits a stepsize that is only inversely proportional to \(\tau\), whereas LT-ADMM-VR requires a stepsize inversely proportional to \(\tau^3\). The latter becomes overly conservative as \(\tau\) grows, leading to slower convergence. This demonstrates that LPDVR is particularly effective when a larger number of local updates is employed.

\subsection{Performance Comparison with Varying Network Size}

\begin{figure}[htpb]
    \begin{minipage}[t]{0.48\linewidth}
        \includegraphics[width=\textwidth]{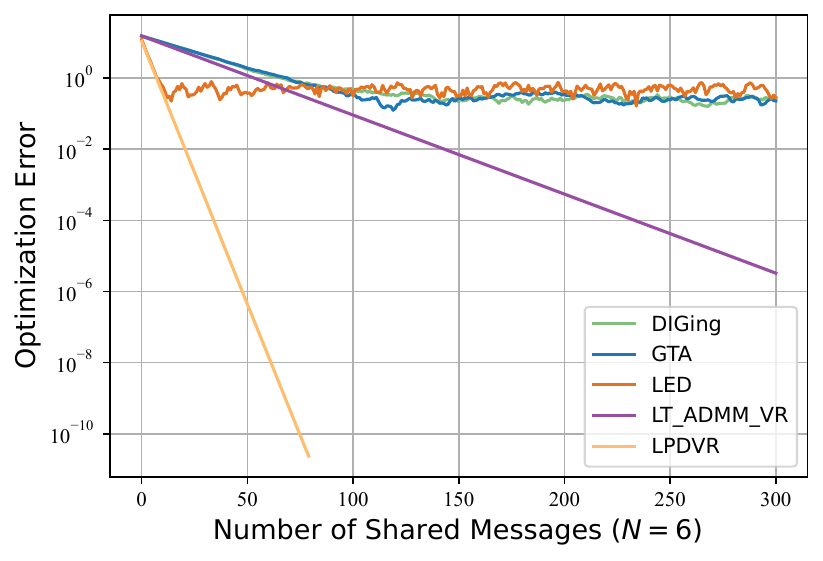}
    \end{minipage}   \hfill
    \begin{minipage}[t]{0.48\linewidth}
        \includegraphics[width=\textwidth]{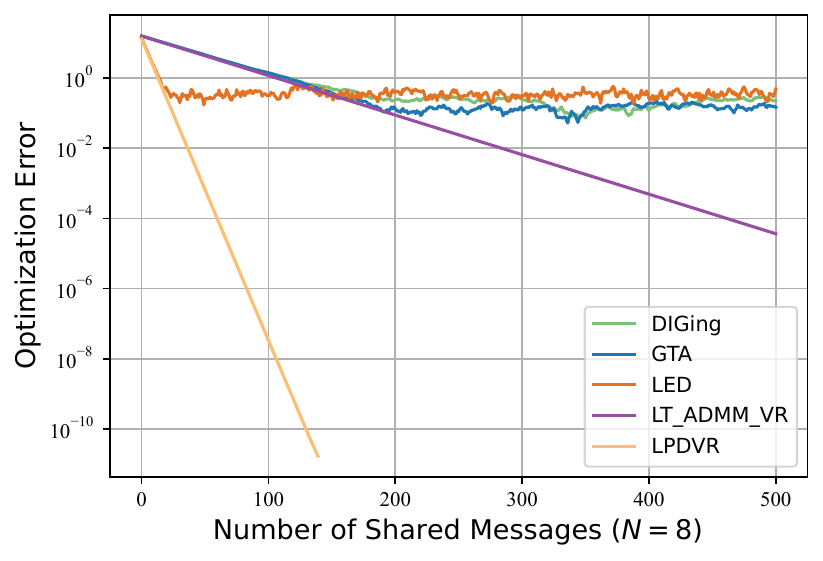}
    \end{minipage}    
    \begin{minipage}[t]{0.48\linewidth}
        \includegraphics[width=\textwidth]{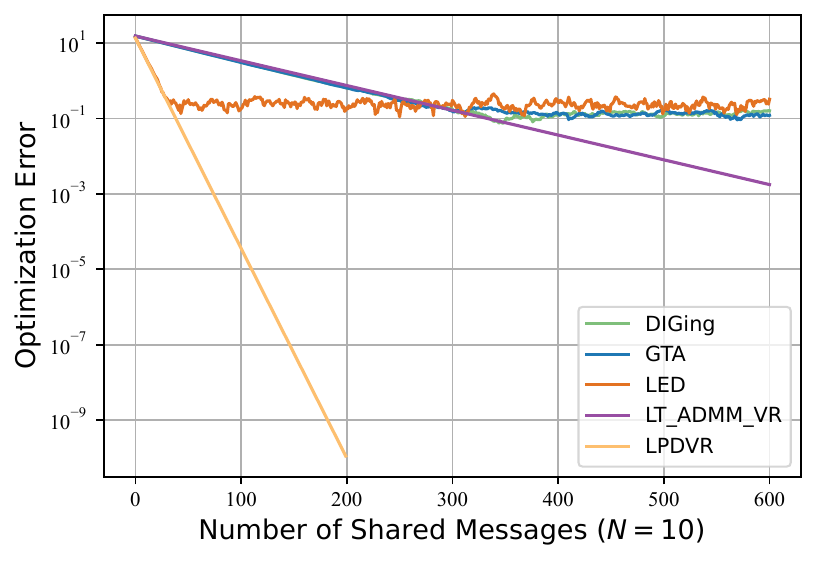}
    \end{minipage}\hfill
    \begin{minipage}[t]{0.48\linewidth}
        \includegraphics[width=\textwidth]{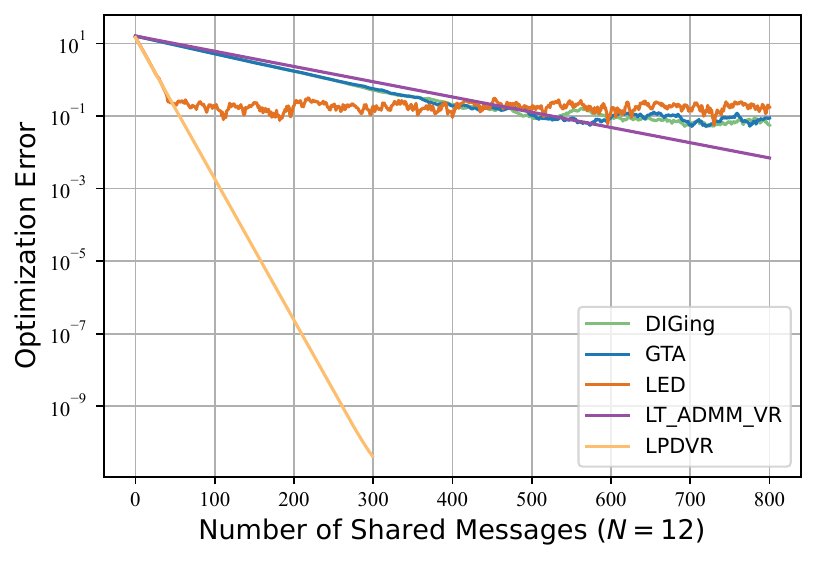}
    \end{minipage}
    \caption{Comparison of LPDVR with LT-ADMM-VR \cite{ren2025communication}, DIGing \cite{WeiShi3}, GTA \cite{LOCAL_SMALL_STEPSIZE}, and LED \cite{LED} over the \textbf{ring} graph under different network sizes $N=6,8,10,12$. }
    \label{ring_different_N}
    \vskip -0.2in
\end{figure}

\begin{figure}[htpb]
    \begin{minipage}[t]{0.48\linewidth}
        \includegraphics[width=\textwidth]{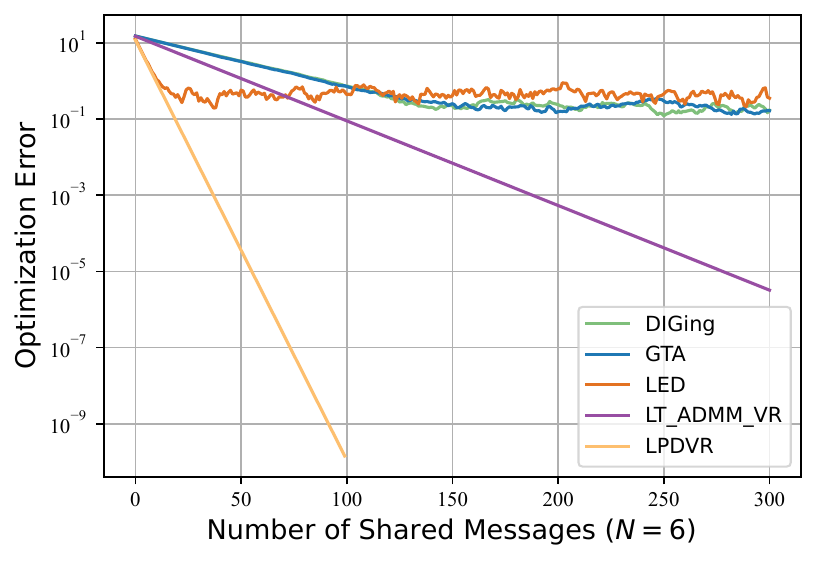}
    \end{minipage}   \hfill
    \begin{minipage}[t]{0.48\linewidth}
        \includegraphics[width=\textwidth]{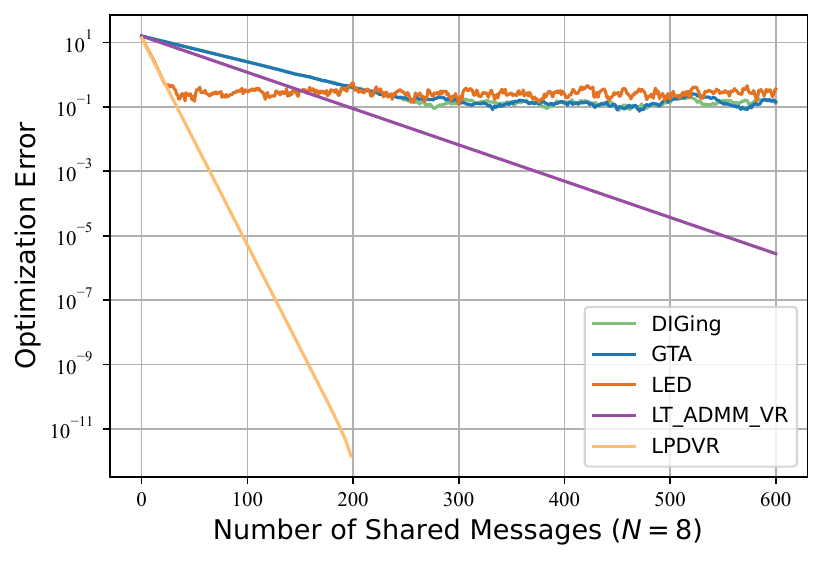}
    \end{minipage}    
    \begin{minipage}[t]{0.48\linewidth}
        \includegraphics[width=\textwidth]{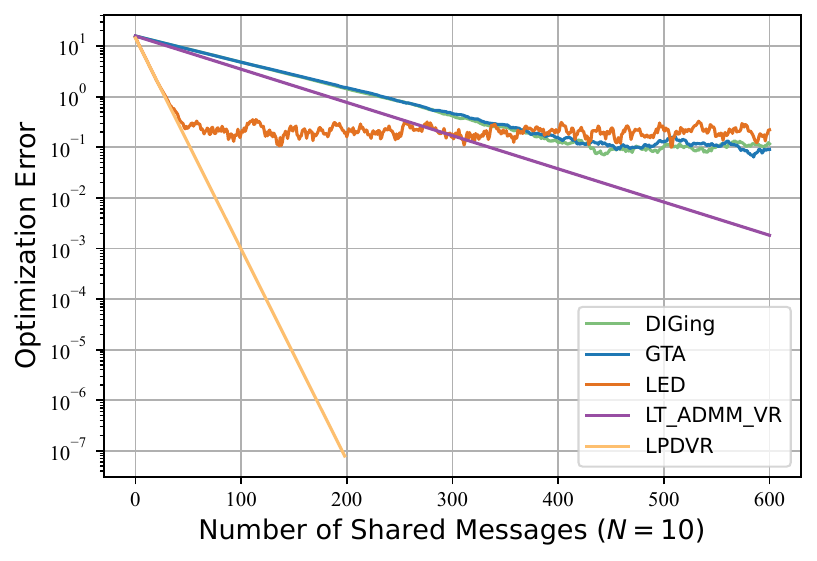}
    \end{minipage}\hfill
    \begin{minipage}[t]{0.48\linewidth}
        \includegraphics[width=\textwidth]{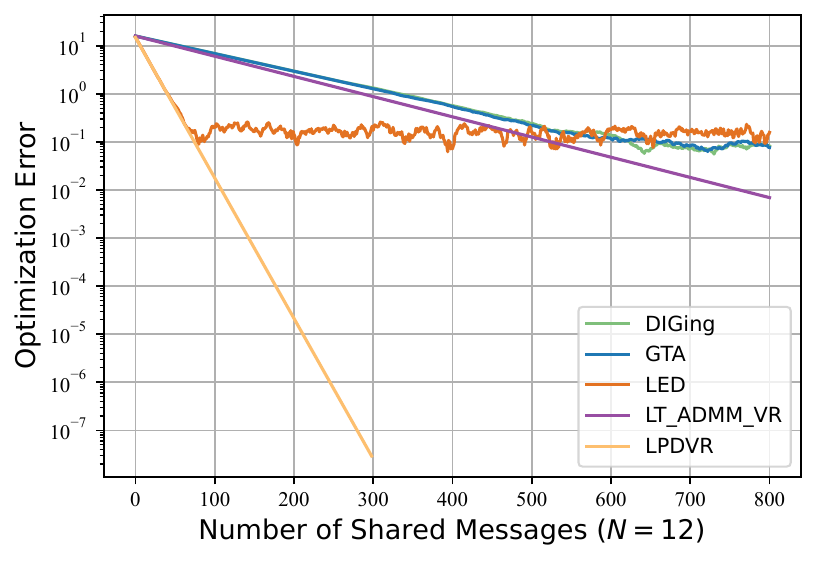}
    \end{minipage}
    \caption{Comparison of LPDVR with LT-ADMM-VR \cite{ren2025communication}, DIGing \cite{WeiShi3}, GTA \cite{LOCAL_SMALL_STEPSIZE}, and LED \cite{LED} over the \textbf{grid} graph under different network sizes $N=6,8,10,12$. }
    \label{grid_different_N}
    \vskip -0.2in
\end{figure}

We further evaluate the scalability of LPDVR by comparing its performance with existing algorithms under different network sizes \(N=6,8,10,12\), with a fixed number of local updates \(\tau=4\). This experiment is also conducted under both ring and grid topologies, and the two settings again lead to consistent observations. As shown in Fig.~\ref{ring_different_N}, under the ring topology, DIGing, GTA, and LED fail to achieve exact convergence across all network sizes, whereas both LPDVR and LT-ADMM-VR consistently converge to high accuracy. Among the two exact-convergent methods, LPDVR consistently outperforms LT-ADMM-VR across all network sizes. Similar observations can be made from Fig.~\ref{grid_different_N} under the grid topology. More importantly, as the number of agents \(N\) increases, the performance gap between LPDVR and LT-ADMM-VR becomes increasingly pronounced under both topologies, with LPDVR maintaining fast convergence while the convergence speed of LT-ADMM-VR deteriorates significantly. This is again consistent with our theoretical results: LPDVR admits a stepsize that is independent of \(N\), whereas LT-ADMM-VR requires a stepsize inversely proportional to \(\sqrt{N}\). The latter becomes overly conservative as the network grows, leading to slower convergence. This demonstrates that LPDVR is well-suited for large-scale networks, a significant advantage in practice where large-scale networks are common.

\section{Conclusion}\label{conclusion}

This paper proposes LPDVR, a novel decentralized optimization algorithm that integrates variance reduction with multiple local updates. We establish $\mathcal{O}(1/T)$ convergence for nonconvex and convex objectives, and linear convergence for strongly convex objectives. This is the first time that exact linear convergence has been achieved for decentralized stochastic optimization with multiple local updates. Beyond convergence, LPDVR offers two practical advantages: its stepsize is independent of the network size and only inversely proportional to the number of local updates, enabling significantly larger stepsizes than existing ADMM-based variance-reduced methods. Moreover, this stepsize condition is comparable to, or even less
restrictive than, those of conventional stochastic gradient or exact gradient methods, demonstrating that exact convergence via variance reduction is achieved without enforcing a more conservative stepsize. In addition, LPDVR requires only one variable exchange per round, halving the per-round communication overhead compared to most existing methods. Numerical results demonstrate the superiority of LPDVR over existing methods in both convergence speed and solution accuracy. Future research may extend LPDVR to nonconvex problems satisfying the PL condition, where linear convergence remains attainable, and to directed graphs, which would enhance its applicability to more realistic network scenarios.

\appendices

\section{Proof of Lemma \ref{lemma_most_important}}\label{lemma_most_important_proof}
Based on the matrix decomposition of $\mathbf{W}$ in Lemma \ref{lemma_from_LED1}, multiplying $\mathbf{Q}^{\top}$ for both sides of (\ref{matrix_form_algorithm_x}) and (\ref{E_update_rule}) yields
\begin{align}\label{appendix_matrix1}
\hspace{-0.3cm}\begin{bmatrix}
\mathbf{Q}^{\top} X(t+1) \\
\mathbf{Q}^{\top} E(t+1)
\end{bmatrix}
=
\begin{bmatrix}
\mathbf{\Lambda} & - \mathbf{\Lambda} \\
\mathbf{B}& \mathbf{\Lambda}
\end{bmatrix}
\begin{bmatrix}
\mathbf{Q}^{\top} {X}(t) \\
\mathbf{Q}^{\top} {E}(t)
\end{bmatrix}
-
\alpha
\begin{bmatrix}
\mathbf{A}_1(t) \\
\mathbf{A}_2(t)
\end{bmatrix},
\end{align}
where $\mathbf{A}_1(t)=\mathbf{\Lambda} \mathbf{Q}^{\top} \sum_{r=0}^{\tau-1} ( \nabla F(\Phi(t,r)) - \nabla F(\bar{\mathbf{X}}(t)) + S(t,r) )$ and $\mathbf{A}_2(t)=\mathbf{B} \mathbf{Q}^{\top} \sum_{r=0}^{\tau-1} ( \nabla F(\Phi(t,r)) - \nabla F(\bar{\mathbf{X}}(t)) + S(t,r) ) -\tau\mathbf{Q}^{\top} ( \nabla F(\bar{\mathbf{X}}(t+1)) - \nabla F(\bar{\mathbf{X}}(t)) )$. Based on Lemma 1 in \cite{LED}, the above matrix equation can be reformulated as 
\begin{align*}
{D}(t+1) 
= \Delta {D}(t)  
- \alpha \mathbf{V}^{-1}
\begin{bmatrix}
\mathbf{A}_3(t) \\
\mathbf{A}_4(t)
\end{bmatrix},
\end{align*}
where $\Delta$ is a matrix with norm $\delta=\Vert \Delta\Vert=\sqrt{\lambda}$, $\mathbf{A}_3(t)=\mathbf{\Lambda} \mathbf{Q}^{\top} \sum_{r=0}^{\tau-1} \bigl(\nabla F(\Phi(t,r)) - \nabla F(\bar{X}(t)) + S(t,r)\bigr)$, and $\mathbf{A}_4(t)=\mathbf{B} \mathbf{Q}^{\top} \sum_{r=0}^{\tau-1} \bigl(\nabla F(\Phi(t,r)) - \nabla F(\bar{X}(t)) + S(t,r)\bigr) - \tau \mathbf{Q}^{\top} \bigl(\nabla F(\bar{X}(t+1)) - \nabla F(\bar{X}(t))\bigr)$. Based on the above equality, Lemma~\ref{lemma_most_important} follows directly from the definitions of $S(t)$ and $H(t)$.

\section{Proof of Lemma \ref{variance_bound_importance1}}\label{variance_bound_importance1_proof}
For $\tau=1$, we have
$\delta(t)=\|\widehat{\Phi}(t)\|^2
=\|X(t)-\bar{X}(t)\|^2$,
so the desired bounds follow directly from
\eqref{D_X_property}.
In the following, assume $\tau\geq2$. With the definition of $\delta(t)$ in \eqref{variance_definition}, we should analyze the properties of $z_{i,h}(t,r)$. From the update rule in \eqref{add_for_if_algorithm}, we obtain
\begin{equation}\label{z_update_rule}
z_{i,h}(t,r+1)=\left\{
\begin{aligned}
\phi_i(t,r+1),\quad &\text{w.p.  }\  \frac{1}{m_i},\\
z_{i,h}(t,r),\quad &\text{w.p. } 1-\frac{1}{m_i}.
\end{aligned}
\right.
\end{equation}
Then, we will analyze the property of $\delta(t)$. We have
\begin{align}
&\mathbb{E}\big[\delta(t,r+1)|\mathcal{F}_{t,r-1}\big]\nonumber\\
%=&\sum^{N}_{i=1}\frac{1}{m_i}\sum^{m_i}_{h=1}\mathbb{E}\big[\Vert z_{i,h}(t,r+1)-\bar{x}(t)\Vert^2|\mathcal{F}_{t,r-1}\big]\nonumber\\
=&\sum^{N}_{i=1}\frac{1}{m_i}\sum^{m_i}_{h=1}\mathbb{E}\big[\frac{1}{m_i}\Vert \phi_i(t,r+1)-\bar{x}(t)\Vert^2\nonumber\\
&+\big(1-\frac{1}{m_i}\big)\Vert z_{i,h}(t,r)-\bar{x}(t)\Vert^2 |\mathcal{F}_{t,r-1}\big]\nonumber\\
\leq&\rho\delta(t,r)+\frac{1}{m}\mathbb{E}\big[\Vert \Phi(t,r+1)-\bar{X}(t)\Vert^2 |\mathcal{F}_{t,r-1}\big],\label{add_more_tac_adding1}
\end{align}
where $\rho=1-\frac{1}{M}$, the first equality follows from \eqref{variance_definition}, the second equality follows by summing
\eqref{add_for_if_algorithm} over the storage indices $h$
and taking the conditional expectation, and the inequality follows from  \eqref{matrix_definition}, (\ref{matrix_definition_average}), \eqref{variance_definition}, together with the relation (\ref{definition_of_Mm}). Taking the total expectation of (\ref{add_more_tac_adding1}) gives
\begin{align}
\mathbb{E}[\delta(t,r+1)]\leq\rho\mathbb{E}[\delta(t,r)]+\frac{1}{m}\mathbb{E}[\Vert \Phi(t,r+1)-\bar{X}(t)\Vert^2 ].\label{delta_step1}
\end{align}
Thus, we should analyze the term $\mathbb{E}[\Vert \Phi(t,r+1)-\bar{X}(t)\Vert^2 ]$ in \eqref{delta_step1}. To this end, we have
\begin{align}
\hspace{-0.1cm}&\mathbb{E}\big[\Vert \phi_i(t,r+1)-\bar{x}(t)\Vert^2\big]\nonumber\\
\hspace{-0.1cm}\leq&\mathbb{E}\big[\Vert \phi_i(t,r)-\alpha  \nabla f_i(\phi_{i}(t,r))-\beta y_i(t)-\bar{x}(t)\Vert^2\big]\nonumber\\
\hspace{-0.1cm}&+\alpha^2\mathbb{E}\big[\Vert \nabla f_i(\phi_{i}(t,r))  -g_i(\phi_{i}(t,r))\Vert^2\big]\nonumber\\
\hspace{-0.1cm}\leq&\frac{\tau}{\tau-1}\mathbb{E}\big[\Vert \phi_i(t,r)-\bar{x}(t)\Vert^2\big]+\tau\mathbb{E}\big[\Vert \alpha  \nabla f_i(\phi_{i}(t,r))\nonumber\\
\hspace{-0.1cm}&+\beta y_i(t)\Vert^2\big]+\alpha^2\mathbb{E}\big[\Vert \nabla f_i(\phi_{i}(t,r))  -g_i(\phi_{i}(t,r))\Vert^2\big],\label{delta_step2}
\end{align}
where the first and second inequalities follow from \eqref{update_phi} and \eqref{expectation_gradient}, and the relation $\Vert x+y\Vert^2\leq (1+a)\Vert x\Vert^2+(1+\frac{1}{a})\Vert y\Vert^2$ with $a=\frac{1}{\tau-1}$, respectively. Furthermore, we have
\begin{align}
&\mathbb{E}\big[\Vert \phi_i(t,r+1)-\bar{x}(t)\Vert^2\big]\nonumber\\
\leq&\frac{\tau}{\tau-1}\mathbb{E}\big[\Vert \phi_i(t,r)-\bar{x}(t)\Vert^2\big]+\tau\mathbb{E}\big[\Vert \alpha  \nabla f_i(\phi_{i}(t,r))+\beta e_i(t)\nonumber\\
&-\alpha  \nabla f_i(\bar{x}(t))\Vert^2\big]+\alpha^2\mathbb{E}\big[\Vert \nabla f_i(\phi_{i}(t,r))  -g_i(\phi_{i}(t,r))\Vert^2\big]\nonumber\\
\leq&\big(\frac{\tau}{\tau-1}+2\tau\alpha^2L^2\big)\mathbb{E}\big[\Vert \phi_{i}(t,r)-\bar{x}(t)\Vert^2\big]+\frac{2}{\tau} \mathbb{E}\big[\Vert e_i(t)\Vert^2\big]\nonumber\\
&+\alpha^2\mathbb{E}\big[\Vert \nabla f_i(\phi_{i}(t,r))  -g_i(\phi_{i}(t,r))\Vert^2\big],\label{delta_step4}
\end{align}
where the first and second inequalities follow from substituting \eqref{matrix_definition2} into \eqref{delta_step2}, and Assumption~\ref{smooth_assumption}, respectively. If the stepsize of LPDVR satisfies $\alpha\leq\frac{1}{2\tau L}$ in Lemma \ref{variance_bound_importance1}, then substituting (\ref{variance_step2}) into \eqref{delta_step4} and summing over $i=1,2,\ldots,N$ yields
\begin{align}
\mathbb{E}\big[\Vert \Phi(t,r&+1)-\bar{X}(t)\Vert^2\big]\leq\frac{\tau+1}{\tau-1}\mathbb{E}\big[\Vert \Phi(t,r)-\bar{X}(t)\Vert^2\big]\nonumber\\
&+\frac{2}{\tau} \mathbb{E}\big[\Vert E(t)\Vert^2\big]+2\alpha^2L^2 \mathbb{E}\big[\delta(t,r)\big].\label{delta_step5_add}
\end{align}
Based on the initial value settings of $\phi_i(t,0)$ during the local update processes in LPDVR, we have $\Phi(t,0)=X(t)$. Furthermore, applying mathematical induction to \eqref{delta_step5_add} yields
\begin{align}
&\mathbb{E}\big[\Vert \Phi(t,r+1)-\bar{X}(t)\Vert^2\big]\leq\xi^{r+1}\mathbb{E}\big[\Vert X(t)-\bar{X}(t)\Vert^2\big]\nonumber\\
&\quad+ \frac{2}{\tau}\mathbb{E}\big[\Vert E(t)\Vert^2\big]\sum^{r}_{l=0}\xi^l+2\alpha^2L^2 \sum^{r}_{l=0}\xi^{r-l}\mathbb{E}\big[\delta(t,l)\big],\label{delta_step6}
\end{align}
where $\xi=1+\frac{2}{\tau-1}$. For $0\leq r\leq\tau-1$ and $0\leq l\leq r$, we have $\max\{\xi^{r+1},\xi^l,\xi^{r-l}\}
\leq \xi^\tau\leq9.$ Combining these bounds with \eqref{delta_step6} yields
\begin{align}
\mathbb{E}\big[\Vert \Phi(t,r&+1)-\bar{X}(t)\Vert^2\big]\leq 9\mathbb{E}\big[\Vert X(t)-\bar{X}(t)\Vert^2\big]\nonumber\\
&+ 18\mathbb{E}\big[\Vert E(t)\Vert^2\big]+18\alpha^2L^2 \sum^{r}_{l=0}\mathbb{E}\big[\delta(t,l)\big].\label{delta_step7}
\end{align}
Next, we will analyze the property of $\sum^{r}_{l=0}\mathbb{E}\big[\delta(t,l)\big]$ in (\ref{delta_step7}). Based on the initial value settings of $z_{i,h}(t,0)$ during the local update processes in LPDVR, we have $\delta(t,0)=\sum^{N}_{i=1}\frac{1}{m_i}\sum^{m_i}_{h=1}\Vert z_{i,h}(t,0)-\bar{x}(t)\Vert^2=\Vert X(t)-\bar{X}(t)\Vert^2$. Based on the above initial value $\delta(t,0)$, applying mathematical induction to \eqref{delta_step1} and then summing over $l=0,1,\cdots,r$ gives $\sum^{r}_{l=0}\mathbb{E}\big[\delta(t,l)\big]\leq\sum^{r}_{l=0}\rho^l \mathbb{E}\big[\Vert X(t)-\bar{X}(t)\Vert^2\big]+\frac{1}{m}\sum^{r}_{l=0}\sum^{l-1}_{s=0}\rho^s\mathbb{E}\big[\Vert \Phi(t,l-s)-\bar{X}(t)\Vert^2 \big].$ Furthermore, applying Lemma \ref{ming_yan_lemma} to this inequality gives
\begin{align}
\sum^{r}_{l=0}\mathbb{E}\big[\delta(t,l)\big]\leq& M\mathbb{E}\big[\Vert X(t)-\bar{X}(t)\Vert^2\big]\nonumber\\
&+\frac{M}{m}\sum^{r}_{l=0}\mathbb{E}\big[\Vert \Phi(t,l)-\bar{X}(t)\Vert^2 \big]\label{delta_step8}
\end{align}
because of $\rho=1-\frac{1}{M}$. Substituting \eqref{delta_step8} into \eqref{delta_step7}
and summing over $r=0,\ldots,\tau-1$, then adding
the initial term
$\mathbb{E}[\|\Phi(t,0)-\bar{X}(t)\|^2]$,
dropping the nonnegative terminal term, and
rearranging terms yields $\big(1-\frac{18\tau\alpha^2L^2M}{m}\big)\mathbb{E}[\Vert \widehat{\Phi}(t)\Vert^2]\leq 36\tau(\mathbb{E}[\Vert E(t)-\bar{E}(t)\Vert^2]+ \mathbb{E}[\Vert \bar{E}(t)\Vert^2])+(1+9\tau+18\tau\alpha^2L^2 M)\mathbb{E}[\Vert X(t)-\bar{X}(t)\Vert^2]$, 
which can be simplified to 
\begin{align}
\mathbb{E}\big[\Vert \widehat{\Phi}(t)\Vert^2\big]\leq &13\tau\mathbb{E}\big[\Vert X(t)-\bar{X}(t)\Vert^2\big]+ 40\tau\mathbb{E}\big[\Vert \bar{E}(t)\Vert^2\big]\nonumber\\
&+ 40\tau\mathbb{E}\big[\Vert E(t)-\bar{E}(t)\Vert^2\big],\label{add_first_3}
\end{align} 
since the stepsize  satisfies $\alpha\leq\min\{\frac{1}{2\tau L},\frac{1}{5\sqrt{M}L}, \frac{1}{30L}\sqrt{\frac{m}{\tau M}}\}$ in Lemma \ref{variance_bound_importance1}. From  Lemma \ref{additional_lemma}, we have $\Vert X(t)-\bar{X}(t)\Vert^2\leq 4 \Vert {D}(t)\Vert^2$ and $\Vert E(t)-\bar{E}(t)\Vert^2\leq 4 \Vert {D}(t)\Vert^2.$ Combining the definition $\bar{E}(t) = \mathbf{1}_N \otimes \bar{e}(t)$ and (\ref{average_e}) gives
\begin{align}
\Vert \bar{E}(t)\Vert^2=\alpha^2\tau^2N \Vert{\nabla f}(\bar{x}(t))\Vert^2.\label{add_first_2}
\end{align}
Substituting (\ref{D_X_property}) and (\ref{add_first_2}) into (\ref{add_first_3}) completes the proof.

\section{Supporting Lemmas}\label{app:supporting_lemmas}
\begin{lemma}[\cite{nesterov2013introductory}]\label{smooth_property_nes}
If $f: \mathbb{R}^n \to \mathbb{R}$ is $L$-smooth over $\mathbb{R}^n$, then $ f(y) \leq f(x) + \langle \nabla f(x), y - x \rangle + \frac{L}{2} \Vert y - x \Vert^2$ holds for all $x, y \in \mathbb{R}^n$.
\end{lemma}
\begin{lemma}[\cite{Ming_Yan2}]\label{ming_yan_lemma}
Let $\{a(t)\}$ and $\{b(t)\}$ be non-negative sequences such that $a(t)=\sum_{s=1}^{t}\rho^{t-s}b(s)$ for some $0\leq \rho<1$. Then the inequality $\sum_{t=1}^{k} a(t) \leq \frac{1}{1-\rho} \sum_{t=1}^{k} b(t)$ holds.
\end{lemma}

\begin{lemma}[\cite{pmlr-v119-karimireddy20a}]\label{scaffold_lemma}
For any \( L \)-smooth and \( \mu \)-strongly convex function \( g \), and arbitrary \( x, y, z \) in its domain, we have
$\langle z-y,\nabla g(x)\rangle \ge g(z)-g(y)+\frac{\mu}{4}\|y-z\|^2 - L\|z-x\|^2.$
\end{lemma}

\begin{lemma}\label{additional_lemma}
The following bounds hold:
\begin{align}
\hspace{-0.23cm}\Vert X(t)-\bar{X}(t) \Vert\leq 2\Vert D(t)\Vert,\ \Vert E(t)-\bar{E}(t) \Vert\leq 2\Vert D(t)\Vert,\label{D_X_property}\\
\hspace{-0.13cm}\mathbb{E}\big[\Vert H(t)\Vert^2\big]\leq 2\alpha^2\tau^2L^2\Vert \mathbf{B}^{-1}\Vert N\big(\mathbb{E}\big[\Vert\sum^{\tau-1}_{r=0}\bar{s}(t,r)\Vert^2\big]\nonumber\\
\hspace{-0.13cm}+\mathbb{E}\big[\Vert\sum^{\tau-1}_{r=0}\overline{\nabla F}(\Phi(t,r))\Vert^2\big]\big)+4\tau L^2\mathbb{E}\big[\Vert\widehat{\Phi}(t)\Vert^2\big],\label{second_part_9}\\
\Vert S(t)\Vert^2\leq 2\Vert \sum^{\tau-1}_{r=0}S(t,r)\Vert^2.\quad\quad\quad\quad\quad\quad\quad\quad\quad\ \label{S_property}
\end{align}
\end{lemma}
\begin{proof}
The proof is similar to that in \cite{LED} and is omitted here for brevity.
\end{proof}

\bibliographystyle{ieeetr}
\bibliography{reference}

\end{document}